\documentclass[11pt,fleqn]{article}
\RequirePackage[OT1]{fontenc}
\RequirePackage{amsthm,amsmath}
\RequirePackage[round]{natbib}
\RequirePackage[colorlinks,citecolor=blue,urlcolor=blue]{hyperref}

\usepackage{amsmath,amsfonts,amssymb,euscript,mathrsfs}
\usepackage{slashbox}
\usepackage{float}
\usepackage[final]{pdfpages}

\usepackage{txfonts}
\usepackage[T1]{fontenc} 
\usepackage{fullpage}
\usepackage{graphicx}
\usepackage{dsfont}
\makeatletter
\def\captionof#1#2{{\def\@captype{#1}#2}}
\makeatother
\def\1{\mbox{\bf 1}}
\def\R{\mathbb{R}}

\def\N{\mathbb{N}}
\def\P{\mathbb{P}}
\def\E{\mathbb{E}}

\def\R{\mathbb{R}}

\def\Z{\mathbb{Z}}

\def\v{\mbox{Var\,}}

\def\c{\mbox{Cov}}

\newtheorem{theo}{Theorem}
\newtheorem{lem}{Lemma}
\newtheorem{prop}{Proposition}
\newtheorem{Def}{Definition}
\newtheorem{cor}{Corollary}
\newtheorem{Def/Prop}{Definition-Proposition}

\newcounter{exos}
\renewcommand\theexos{\arabic{exos}}

\newcounter{prob}
\renewcommand\theprob{\arabic{prob}}

\begin{document}
\author{Lionel Truquet \footnote{UMR 6625 CNRS IRMAR, University of Rennes 1, Campus de Beaulieu, F-35042 Rennes Cedex, France and}
\footnote{Campus de Ker-Lann, rue Blaise Pascal, BP 37203, 35172 Bruz cedex, France. {\it Email: lionel.truquet@ensai.fr}.}
 }
\title{Local stationarity and time-inhomogeneous Markov chains}
\date{}
\maketitle

\begin{abstract}
In this paper, we study a notion of local stationarity for discrete time Markov chains which is useful 
for applications in statistics. In the spirit of some locally stationary processes introduced in the literature, we consider triangular arrays of time-inhomogeneous Markov chains, defined by some families of contracting Markov kernels.
Using the Dobrushin's contraction coefficients for various metrics, we show that the distribution of such Markov chains can be approximated locally with the distribution of ergodic Markov chains and we also study some mixing properties. 
From our approximation results in Wasserstein metrics, we recover several properties obtained for autoregressive processes. Moreover, using the total variation distance or more generally some distances induced by a drift function, we consider new models, such as finite state space Markov chains with time-varying transition matrices or some time-varying versions of integer-valued autoregressive processes. For these two examples, nonparametric kernel estimation of the transition matrix is discussed.
\end{abstract}

\section{Introduction}
Time-inhomogeneous Markov chains have received much less attention in the literature than the homogeneous case. Such chains have been studied mainly for their long-time behavior, often in connexion with the convergence of stochastic algorithms. An introduction to inhomogeneous Markov chains and their use in Monte Carlo methods can be found in \citet{Winkler}. More recent quantitative results for their long time behavior can be found for instance in \citet{Mouli}, \citet{SC1}, or \citet{SC2}. 
In this paper, we consider convergence properties of nonhomogeneous Markov chains but with a different perspective, motivated by applications in mathematical statistics and in the spirit of the notion of local stationarity introduced by \citet{Dahlhaus}. Locally stationary processes have received a considerable attention over the last twenty years, in particular for their ability to model data sets for which time-homogeneity is unrealistic. Locally stationary autoregressive processes (here with one lag for simplicity) can be defined by modifying a recursive equation followed by a stationary process. If $(X_k)_{k\in\Z}$ is a stationary processes defined by $X_k=F_{\theta}\left(X_{k-1},\varepsilon_k\right)$, where $(\varepsilon_k)_{k\in\Z}$ is a sequence of i.i.d random variables and $\theta\in \Theta$ is a parameter, its locally stationary version is usually defined recursively by 
$$X_{n,k}=F_{\theta(k/n)}\left(X_{n,k-1},\varepsilon_k\right),\quad 1\leq k\leq n,$$
where $\theta:[0,1]\rightarrow \Theta$ is a smooth function. This formalism was exploited for defining locally stationary versions of classical time-homogeneous autoregressive processes. See for instance \citet{DR}, \citet{SR} or \citet{Vogt}.  
The term local stationarity comes from the fact that, under some regularity conditions, if $k/n$ is close to a point $u$ of $[0,1]$, 
$X_{n,k}$ is close in some sense to $X_k(u)$ where $\left(X_k(u)\right)_{k\in \Z}$ is the stationary process defined by 
$$X_k(u)=F_{\theta(u)}\left(X_{k-1}(u),\varepsilon_k\right),\quad k\in \Z.$$
Though local stationary processes defined recursively are examples of time-inhomogeneous Markov chains, the properties of these processes are usually derived using this particular autoregressive representation and without exploiting the link with Markov chains. This is one the main difference with respect to stationary processes for which the connection between autoregressive processes and Markov chains has been widely used. See for example the classical textbook of \citet{MT} for many examples of iterative systems studied using Markov chains properties. As a limitation, the simple case of a locally stationary version of finite state space Markov chains has not been considered in the literature.

In this paper, we consider general Markov chains models which will generalize the existing (Markovian) locally stationary processes. Since we do not work directly with autoregressive representations, our definition of local stationarity is based on the approximation of the finite dimensional distributions of the chain with that of some ergodic Markov chains. Let us now give the framework used in the rest of the paper.   
Let $\left(E,d\right)$ be a metric space, $\mathcal{B}(E)$ its corresponding Borel $\sigma-$field and $\left\{Q_u: u\in[0,1]\right\}$  a family of Markov kernels on $\left(E,\mathcal{B}(E)\right)$. By convention, we set $Q_u=Q_0$ when $u<0$.
We will consider triangular arrays $\left\{X_{n,j}: j\leq n,n\in\Z^+\right\}$ such that for all $n\in\Z^+$, the sequence $\left(X_{n,j}\right)_{j\leq n}$ is a non homogeneous Markov chain such that
$$\P\left(X_{n,k}\in A\vert X_{n,k-1}=x)\right)=Q_{k/n}(x,A),\quad k\leq n .$$
In the sequel the family $\left\{Q_u: u\in [0,1]\right\}$ of Markov kernels will always satisfy some regularity conditions and contraction properties. Precise assumptions will be given in three following sections, but from now on, we assume here that for all $u\in [0,1]$, $Q_u$ has a single invariant probability denoted by $\pi_u$.  
For all positive integer $j$ and all integer $k$ such that $k+j-1\leq n$, we denote by $\pi^{(n)}_{k,j}$ the probability distribution of the vector $\left(X_{n,k},X_{n,k+1},\ldots,X_{n,k+j-1}\right)$ and by $\pi_{u,j}$ the corresponding finite dimensional distribution for the ergodic chain with Markov kernels $Q_u$.  
Loosely speaking, the triangular array will be said locally stationary if for all positive integer $j$, the probability distribution $\pi^{(n)}_{k,j}$ is close to $\pi_{u,j}$ when the ratio $k/n$ is close to $u$. 
A formal definition is given below. For an integer $j\geq 1$, we denote by $\mathcal{P}(E^j)$ the set of probability measures on $E^j$.
\begin{Def}\label{debase}
The triangular array of non-homogeneous Markov chains $\left\{X_{n,k}, n\in\Z^+, k\leq n\right\}$ is said to be locally stationary if for all integer $j\geq 1$, there exists a metric $\vartheta_j$ on 
$\mathcal{P}\left(E^j\right)$, metrizing the topology of weak convergence, such that the two following conditions are satisfied.
\begin{enumerate}
\item
The application $u\mapsto \pi_{u,j}$ is continuous.
\item
$\displaystyle\lim_{n\rightarrow\infty}\sup_{k\leq n-j+1}\vartheta_j\left(\pi_{k,j}^{(n)},\pi_{\frac{k}{n},j}\right)=0$
\end{enumerate}
\end{Def}
In particular, under the two conditions of Definition \ref{debase}, for all continuous and bounded function $f:E^j\rightarrow \R$ and some integers $k=k_n\leq n-j+1$ such that $\lim_{n\rightarrow \infty} k/n =u\in[0,1]$, we have 
$$\lim_{n\rightarrow \infty}\E f\left(X_{n,k},\ldots,X_{n,k+j-1}\right)=\lim_{n\rightarrow \infty}f d\pi^{(n)}_{k,j}=\E f\left(X_1(u),\ldots,X_j(u)\right)=\int f d\pi_{u,j},$$
where $\left(X_k(u)\right)_{k\in\Z}$ denotes a stationary Markov chain with transition $Q_u$.
In this paper, Condition $1$ will always hold from the Hölder continuity properties that we will assume for the application $u\mapsto Q_u$.
Of course, the metrics $\vartheta_j$ will of the same nature for different integers $j$, e.g the total variation distance on $\mathcal{P}(E^j)$.

In this paper, we will consider three type of metrics on $\mathcal{P}(E)$ for approximating $\pi^{(n)}_k$ by $\pi_{k/n}$ or $\pi_u$ (and in a second step for approximating an arbitrary finite dimensional distribution) and deriving mixing properties of these triangular arrays.
We will extensively make use of the so-called Dobrushin's contraction coefficient.
 In Section \ref{Dob1}, we consider the total variation distance. 
This is the metric for which the contraction coefficient for Markov kernels has been originally introduced by \citet{Dob}.
Contraction properties of the kernels $Q_u$ or their iteration with respect to this metric will enable us to consider a model
of nonhomogeneous finite state space Markov chains for which we will study a nonparametric estimator of the time-varying transition matrix. In Section \ref{Dob2}, we consider contraction properties for Wasserstein metrics. The contraction coefficient for the Wasserstein metric of order $1$ has been first considered by \citet{Dob+} for giving sufficient conditions under which a system of conditional distributions defines a unique joint distribution. We will consider more generally the Wasserstein metric of order $p\geq 1$. 
This type of metric is very well adapted for recovering some results obtained for autoregressive processes with time-varying coefficients. Finally, in Section \ref{Dob3}, we consider Markov kernels 
satisfying drift and minoration conditions ensuring geometric ergodicity and for which \citet{HM1} have recently found a contraction property for a metric induced by a modified drift function. We illustrate this third approach with the statistical inference of some integer-valued autoregressive processes with time-varying coefficients.

\section{Total variation distance and finite state space Markov chains}\label{Dob1}
Let us first give some notations that we will extensively use in the sequel. If $\mu\in\mathcal{P}(E)$ and $R$ is a probability kernel from $\left(E,\mathcal{B}(E)\right)$ to $\left(E,\mathcal{B}(E)\right)$, we will denote by $\mu R$ the probability measure defined by 
$$\mu R(A)=\int R(x,A)d\mu(x),\quad A\in\mathcal{B}(E).$$
Moreover if $f:E\rightarrow \R$ is a measurable function, we set 
$\mu f=\int fd\mu$ and $Rf:E\rightarrow \R$ will be the function defined by
$Rf(x)=\int R(x,dy)f(y)$, $x\in E$, provided these integrals are well defined.
Finally, the Dirac measure at point $x\in E$ is denoted by $\delta_x$.
\subsection{Contraction and approximation result for the total variation distance}
The total variation distance between two probability measures $\mu,\nu\in\mathcal{P}(E)$ is defined by  
$$\Vert \mu-\nu\Vert_{TV}=\sup_{A\in \mathcal{B}(E)}\left\vert \mu(A)-\nu(A)\right\vert=\frac{1}{2}\sup_{\Vert f\Vert_{\infty}\leq 1}\left\vert \int fd\mu-\int fd\nu\right\vert,$$
where for a measurable function $f:E\rightarrow \R$,  $\Vert f\Vert_{\infty}=\sup_{x\in E}\left\vert f(x)\right\vert$.
  
For the family $\left\{Q_u: u\in [0,1]\right\}$, the following assumptions will be needed.
\begin{description}
\item [A1]
There exist an integer $m\geq 1$ and $r\in (0,1)$ such that for all $(u,x,y)\in [0,1]\times E^2$, 
$$\Vert\delta_x Q_u^m-\delta_y Q_u^m\Vert_{TV}\leq r.$$

\item [A2]
There exist a positive real number $L$ and $\kappa\in (0,1)$ such that for all $(u,v,x)\in [0,1]^2\times E$, 
$$\Vert\delta_x Q_u-\delta_x Q_v\Vert_{TV}\leq L \vert u-v\vert^{\kappa}.$$
\end{description}

The Dobrushin contraction coefficient of Markov kernel $R$ on $\left(E,\mathcal{B}(E)\right)$ is defined by \\$c(R)=\sup_{(x,y)\in E^2}\Vert \delta_x R-\delta_y R\Vert_{TV}$. We have $c(R)\in [0,1]$. Hence, assumption ${\bf A1}$ means that \\$\sup_{u\in [0,1]}c\left(Q_u^m\right)<1$. 
We will still denote by $\Vert \cdot\Vert_{TV}$ the total variation distance (or the total variation norm if we consider the space of signed measures) on $\mathcal{P}(E^j)$ for any integer $j$.
Moreover, let $\left(X_k(u)\right)_{k\in \Z}$ be a stationary Markov chain with transition $Q_u$, for $u\in [0,1]$.
We remind that for an integer $j\geq 1$, $\pi^{(n)}_{k,j}$ (resp. $\pi_{u,j}$) denotes the probability distribution of the vector $\left(X_{n,k},\ldots,X_{n,k+j-1}\right)$ (resp. of the vector $\left(X_k(u),\ldots,X_{k+j-1}(u)\right)$),
\begin{theo}\label{first}
Assume that assumptions ${\bf A1-A2}$ hold true. Then for all $u\in [0,1]$, the Markov kernel $Q_u$ has a single invariant probability $\pi_u$. The triangular array of Markov chain $\left\{X_{n,k}, n\in\Z^+,k\leq n\right\}$ is locally stationary.
Moreover, there exists a positive real number $C$, only depending on $L,m,r,\kappa$ such that 
$$\Vert\pi^{(n)}_{k,j}-\pi_{u,j}\Vert_{TV}\leq C\left[\sum_{s=k}^{k+j-1}\left\vert u-\frac{s}{n}\right\vert^{\kappa}+\frac{1}{n^{\kappa}}\right].$$
\end{theo}
 \paragraph{Note.} Assumption ${\bf A1}$ is satisfied if there exist a positive real number $\varepsilon$, a positive integer $m$ and a family of probability measures $\left\{\nu_u:u\in [0,1]\right\}$ such that  
$$Q^m_u(x,A)\geq \varepsilon \nu_u(A),\mbox{ for all } (u,x,A)\in [0,1]\times E\times \mathcal{B}(E).$$
In the homogeneous case, this condition is the so-called Doeblin's condition (see \citet{MT}, Chapter $16$ for a discussion about this condition). To show that this condition is sufficient for ${\bf A1}$, one can use the inequalities
$$Q_u^m(x,A)-Q_u^m(y,A)\leq 1-\varepsilon +\varepsilon\nu_u(E\setminus A)-Q^m_u(x,E\setminus A)\leq 1-\varepsilon.$$
For a Markov chain with a finite state space, the Doeblin's condition is satisfied if 
$\inf_{u\in [0,1]} Q_u^m(x,y)>0$, taking the counting measure for $\nu_u$. More generally, this condition is satisfied if $Q_u(x,A)=\int_A f_u(x,y)\nu(dy)$ with a probability measure $\nu$ and a density uniformly lower bounded, i.e $\varepsilon=\displaystyle\inf_{(u,x,y)\in [0,1]\times E^2}f_u(x,y)>0$.

\paragraph{Proof of Theorem \ref{first}}
We remind that for a Markov kernel $R$ on $\left(E,\mathcal{E}\right)$ and $\mu,\nu\in\mathcal{P}(E)$, we have 
$$\Vert \mu R-\nu R\Vert_{TV}\leq c(R)\cdot\Vert\mu-\nu\Vert_{TV},$$
where $c(R)=\sup_{(x,y)\in E}\Vert \delta_x R-\delta_y R\Vert_{TV}\in [0,1]$.
Then, under our assumptions, the application $T:\mathcal{P}(E)\rightarrow \mathcal{P}(E)$ defined by $T(\mu)=\mu Q_u^m$ is contractant and
the existence and uniqueness of an invariant probability $\pi_u$ easily follows from the fixed point theorem in a complete metric space. 

We next show Condition $1$ of Definition \ref{debase}. The result is shown by induction. For $j=1$, we have from assumption ${\bf A1}$,
\begin{eqnarray*}
\Vert \pi_u-\pi_v\Vert_{TV}&\leq &\Vert \pi_uQ^m_u-\pi_v Q^m_u\Vert_{TV}+\Vert \pi_v Q^m_u-\pi_v Q^m_v\Vert_{TV}\\
&\leq& r\Vert \pi_u-\pi_v\Vert_{TV}+\sup_{x\in E}\Vert \delta_x Q_u^m-\delta_x Q_v^m\Vert_{TV}.
\end{eqnarray*}
Since for two Markov kernels $R$ and $\widetilde{R}$ and $\mu,\nu\in\mathcal{P}(E)$, we have 
$$\Vert \mu R-\nu\widetilde{R}\Vert_{TV}\leq \sup_{x\in E}\Vert \delta_xR-\delta_x\widetilde{R}\Vert_{TV}+ c\left(\widetilde{R}\right)\Vert \mu-\nu\Vert_{TV},$$ 
we deduce from assumption ${\bf A2}$ that $\sup_{x\in E}\Vert \delta_x Q_u^m-\delta_x Q_v^m\Vert_{TV}\leq mL\vert u-v\vert^{\kappa}$. 
This leads to the inequality
$$\Vert \pi_u-\pi_v\Vert_{TV}\leq \frac{m L}{1-r}\vert u-v\vert^{\kappa}$$
which shows the result for $j=1$. If the continuity condition holds true for $j-1$, we note that $$\pi_{u,j}\left(dx_1,\ldots,d_{x_{j-1}}\right)=\pi_{u,j-1}\left(dx_1,\ldots,d_{x_{j-1}}\right)Q_u\left(x_{j-1},d x_j\right).$$
Moreover, we have 
$$\Vert\pi_{u,j}-\pi_{v,j}\Vert_{TV}\leq \sup_{x\in E}\Vert \delta_x Q_u-\delta_x Q_v\Vert_{TV}+ \Vert \pi_{u,j-1}-\pi_{v,j-1}\Vert_{TV},$$
which leads to the continuity of $u\mapsto \pi_{u,j}$. This justifies Condition $1$ of Definition \ref{debase}.

Finally we prove the bound announced for $\Vert \pi^{(n)}_{k,j}-\pi_{u,j}\Vert_{TV}$.
Note that this bound automatically implies Condition $2$ of Definition \ref{debase}. Let us first note that if
$R_{k,m}=Q_{\frac{k-m+1}{n}}Q_{\frac{k-m+2}{n}}\cdots Q_{\frac{k}{n}}$, we have from assumption ${\bf A2}$, 
$$\sup_{x\in E}\Vert \delta_x R_{k,m}-\delta_x Q_u^m\Vert_{TV}\leq L\sum_{s=k-m+1}^k\left\vert u-\frac{s}{n}\right\vert^{\kappa}.$$ 
Now for $j=1$, we have 
\begin{eqnarray*}
\Vert \pi^{(n)}_k-\pi_u\Vert_{TV}&\leq& \Vert \pi^{(n)}_{k-m}R_{k,m}-\pi^{(n)}_{k-m}Q_u^m\Vert_{TV}+\Vert \pi^{(n)}_{k-m} Q_u^m-\pi_u Q_u^m\Vert_{TV}\\
&\leq& L\sum_{s=k-m+1}^k\left\vert u-\frac{s}{n}\right\vert^{\kappa}+r \Vert \pi^{(n)}_{k-m}-\pi_u\Vert_{TV}.
\end{eqnarray*}
Using the fact that is $s\leq 0$, $\vert u-s/n\vert\leq \vert u\vert$, we deduce that 
$$\Vert \pi^{(n)}_k-\pi_u\Vert_{TV}\leq L\sum_{\ell=0}^{\infty} r^{\ell}\sum_{s=k-(\ell+1)m+1}^{k-\ell m}\left\vert u-\frac{s}{n}\right\vert^{\kappa},$$
which shows the result for $j=1$. Next, using the same argument as for the continuity of the finite-dimensional distributions, we have 
$$\Vert \pi^{(n)}_{k,j}-\pi_{u,j}\Vert_{TV}\leq L\left\vert u-\frac{k+j-1}{n}\right\vert^{\kappa}+\Vert \pi^{(n)}_{k,j-1}-\pi_{u,j-1}\Vert_{TV}.$$
Hence the result easily follows by iteration.$\square$ 

\subsection{$\beta-$mixing properties}
In this subsection, we consider the problem of mixing for the locally stationary Markov chains introduced previously. 
For convenience, we assume that $X_{n,j}$ is equal to zero if $j\geq n+1$. 
For a positive integer $n$ and an integer $i\in\Z$, we denote by $\mathcal{F}^{(n)}_i$ the sigma field $\sigma\left(X_{n,j}: j\leq i\right)$. Now setting 
$$V\left(\mathcal{F}^{(n)}_i,X_{n,i+j}\right)=\sup\left\{\left\vert \E\left[f\left(X_{n,i+j}\vert \mathcal{F}^{(n)}_i\right)\right]-\E\left[f\left(X_{n,i+j}\right)\right]\right\vert:\quad f\mbox{ s.t } \Vert f\Vert_{\infty}\leq 1\right\},$$
the $\beta_n-$mixing coefficient for the sequence $\left(X_{n,j}\right)_{j\in\Z}$ 
is defined by 
$$\beta_n(j)=\frac{1}{2}\sup_{i\in\Z}\E\left[V\left(\mathcal{F}^{(n)}_i,X_{n,i+j}\right)\right].$$
Under our assumptions, this coefficient is shown to decrease exponentially fast.

\begin{prop}\label{mixing1}
Assume that assumptions ${\bf A1-A2}$ hold true. Then there exist $C>0$ and $\rho\in (0,1)$, only depending on $m,L,\kappa$ and $r$ such that 
$$\beta_n(j)\leq C \rho^{[j/m]},$$
where $[x]$ denotes the integer part of a real number $x$.
\end{prop}

\paragraph{Note.} The usual strong mixing coefficient is defined for Markov chains by 
$$\alpha_n(j)=\sup_{i\in\Z}\left\{\left\vert \P(A\cap B)-\P(A)\P(B)\right\vert: A\in \sigma\left(X_{n,i}\right), B\in\sigma\left(X_{n,i+j}\right)\right\}.$$
We have $\alpha_n(j)\leq \beta_n(j)$. We refer the reader \citet{Do} for the definition of some classical mixing coefficients and their properties. In this paper, we will mainly use some results available for the larger class of strong-mixing processes.

\paragraph{Proof of Proposition \ref{mixing1}}
We first consider $\epsilon>0$ such that $\rho=2m L\epsilon^{\kappa}+r<1$. Assume first that $n\geq \frac{m}{\epsilon}$. 
For $k\leq n$, we set $Q_{k,m}=Q_{\frac{k-m+1}{n}}\cdots Q_{\frac{k}{n}}$. By noticing that under Assumption ${\bf A2}$, we have
$$\sup_{\mu\in\mathcal{P}(E)}\Vert\mu Q_u-\mu Q_v\Vert_{TV}\leq L\vert u-v\vert^{\kappa},$$
we deduce the bound 
$$\sup_{x\in E}\Vert \delta_x Q_{k,m}-\delta_x Q^m_{\frac{k}{n}}\Vert_{TV}\leq  m L\epsilon^{\kappa}.$$
Then, from Assumption ${\bf A1}$, we get 
$$\sup_{x,y\in E}\Vert \delta_x Q_{k,m}-\delta_y Q_{k,m}\Vert_{TV}\leq \rho.$$
Now if $j=tm+s$ for two positive integers $t,s$, we get 
$$\Vert \delta_{X_{n,k-j}}Q_{\frac{k-j+1}{n}}\cdots Q_{\frac{k}{n}}-\pi_{k-j}^{(n)}Q_{\frac{k-j+1}{n}}\cdots Q_{\frac{k}{n}}\Vert_{TV}\\
\leq \rho^t.$$
Now, if $n<\frac{m}{\epsilon}$, one can show that $\beta_n(j)\leq 1$, if $j\leq n$ and $\beta_n(j)\leq r^{\left[\frac{j-n}{m}\right]}$ if $j>n$. This leads to the result with an appropriate choice of $C$, e.g $C=\rho^{-\frac{1}{\epsilon}-1}$.$\square$

\subsection{Finite state space Markov chains}\label{poid}
Let $E$ be a finite set. In this case, we obtain the following result.
\begin{cor}\label{finite}
Let $\left\{Q_u: u\in[0,1]\right\}$ be a family of transition matrices such that for each $u\in [0,1]$, the Markov chain with transition matrix $Q_u$ is irreducible and aperiodic. 
Assume further that for all $(x,y)\in E^2$, the application $u\rightarrow Q_u(x,y)$ is $\kappa-$Hölder continuous. Then Theorem \ref{central} applies and the $\beta-$mixing coefficients are bounded as in Proposition \ref{mixing1}.$\square$ 
\end{cor}

\paragraph{Proof of Corollary\ref{finite}}
Using the fact that 
$$\Vert\delta_x Q_u^m-\delta_y Q_u^m\Vert_{TV}=1-\sum_{z\in E}Q^m_u(x,z)\wedge Q_u^m(y,z)\leq 1-\vert E\vert\cdot\inf_{x,y\in E}Q_u^m(x,y).$$
Then assumption ${\bf A1}$ is satisfied as soon as $\inf_{u\in[0,1],(x,y)\in E^2}Q^m_u(x,y)>0$. From aperiodicity and irreducibility, it is well know that for each $u\in [0,1]$, 
$$m_u=\inf\left\{k\geq 1: \min_{(x,y)\in E^2}Q_u^{k}(x,y)>0\right\}<\infty.$$
By continuity, the sets $\mathcal{O}_u=\left\{v\in [0,1]: P_v^{m_u}>0\right\}$ are open subsets of $[0,1]$. From the compactness of the interval $[0,1]$, 
$[0,1]$ can be covered by finitely many $\mathcal{O}_u$, say $\mathcal{O}_{u_1},\ldots,\mathcal{O}_{u_d}$. Then assumption ${\bf A1}$ is satisfied with $m=\max_{1\leq i\leq d}m_{u_i}$. Assumption ${\bf A2}$ is automatically satisfied and Theorem \ref{central} applies.$\square$

Now, we show that our results can be used for nonparametric kernel estimation of the invariant probability $\pi_u$ or the transition matrix $Q_u$. This kind of estimation requires an estimation of quantities of type $h_u=\E\left[f\left(X_1(u),\ldots,X_{\ell}(u)\right)\right]$ where $f:E^{\ell}\rightarrow \R$ is a function and $\ell$ is an integer. To this end, a classical method used for locally stationary time series is based on kernel estimation. 
See for instance \citet{DR}, \citet{Fryz}, \citet{Vogt} or \citet{ZW} for nonparametric kernel estimation of locally stationary processes. Let  
$K:\R\rightarrow \R_+$ be a Lipschitz function supported on $[-1,1]$ and such that $\int K(z)dz=1$. For $b=b_n\in (0,1)$, we set 
$$e_i(u)=\frac{\frac{1}{nb}K\left(u-\frac{i}{n}\right)}{\frac{1}{nb}\sum_{j=\ell}^n K\left(u-\frac{j}{n}\right)},\quad u\in [0,1],\quad \ell\leq i\leq n.$$
A natural estimator of $h_u$ is 
$$\hat{h}_u=\sum_{i=\ell}^ne_i(u)f\left(X_{n,i-\ell+1},\ldots,X_{n,i}\right).$$
The next proposition gives a uniform control of the variance part $\hat{h}_u-\E\hat{h}_u$. 
\begin{prop}\label{exp}
Assume that assumption {\bf A3} holds true and that $b\rightarrow 0$, $nb^{1+\epsilon}\rightarrow \infty$ for some $\epsilon>0$. Then
$$\sup_{u\in [0,1]}\left\vert \hat{h}_u-\E \hat{h}_u\right\vert=O_{\P}\left(\frac{\sqrt{\log n}}{\sqrt{nb}}\right).$$ 
\end{prop}

\paragraph{Proof of Proposition \ref{exp}}
We set $Y_{n,i}=f\left(X_{n,i-\ell+1},\ldots,X_{n,i}\right)$. 
First, note that the triangular array $\left(Y_{n,i}\right)_{1\leq i\leq n}$ is $\beta-$mixing (and then $\alpha-$mixing)
with $\beta_n(j)\leq \widetilde{C} \rho^{\left[\frac{j-\ell}{m}\right]}\leq \widetilde{C}\rho^{-1-\frac{\ell}{m}}\widetilde{\rho}^j$ where $\widetilde{C}$ is a positive constant and $\widetilde{\rho}=\rho^{1/m}$. 
We have $[0,1]=\cup_{s=1}^{k+1}I_s$ where $k$ is the integer part of $1/b$, $I_s=((s-1)b,sb]$ for $1\leq s\leq k$ and $I_{k+1}=(kb,1]$.
We set $S_0^{(n)}=0$ and if $1\leq i\leq n$, $S_i^{(n)}=\sum_{s=1}^i Z_s^{(n)}$, where $Z_s^{(n)}=Y_{n,s}-\E Y_{n,s}$. 
Then for $1\leq j\leq j+k\leq n$, we have 
\begin{eqnarray*}
\left\vert \sum_{i=j}^{j+k}e_i(u)Z_i^{(n)}\right\vert &\leq& e_j(u)\cdot\left\vert S_{j-1}^{(n)}\right\vert+e_{j+k}(u)\cdot\left\vert S_{j+k}^{(n)}\right\vert+\sum_{i=j}^{j+k-1}\left\vert e_i(u)-e_{i-1}(u)\right\vert\cdot \left\vert S_i^{(n)}\right\vert\\
&\leq & \frac{C''}{nb}\max_{j-1\leq i\leq j+k}\left\vert S_i^{(n)}\right\vert.
\end{eqnarray*}
This gives the bound 
\begin{eqnarray*}
\max_{u\in [0,1]}\left\vert \sum_{i=1}^n e_i(u) Z_i^{(n)}\right\vert&\leq& \max_{1\leq s\leq k}\max_{u\in I_s}\left\vert \sum_{n(s-2)b\leq i\leq n(s+1)b}e_i(u)Z_i^{(n)}\right\vert\\
&\leq& \frac{C''}{nb}\max_{1\leq s\leq k+1} \max_{n(s-2)b-1\leq i\leq n(s+1)b}\left\vert S_i^{(n)}\right\vert.
\end{eqnarray*}
We will use the exponential inequality for strong mixing sequences given in \citet{Rio1}, Theorem $6.1$ (see also \citet{Rio11}, Theorem $6.1$). This inequality guarantees that for any integer $q$, we have 
$$\P\left(\max_{n(s-2)b-1\leq i\leq n(s+1)b}\left\vert S_i^{(n)}\right\vert\geq F\lambda\right)\leq G\exp\left(-\frac{\lambda}{2q\Vert f\Vert_{\infty}}\log\left(1+K\frac{\lambda q}{nb}\right)\right)+Mnb\frac{\widetilde{\rho}^q}{\lambda},$$
where $F,G,K,M$ are three positive real numbers not depending on $n$ and $s$ and $\lambda\geq q \Vert f\Vert_{\infty}$.
We have $k=O\left(b^{-1}\right)$ and setting $q\approx \frac{\sqrt{nb}}{\sqrt{\log n}}$ and $\lambda=\lambda' \sqrt{nb \log n}$, we have for $\lambda'$ large enough
$$\P\left(\max_{u\in [0,1]}\left\vert \sum_{i=1}^n e_i(u) Z_i^{(n)}\right\vert>\frac{F \lambda}{nb}\right)=O_{\P}\left(\frac{1}{bn^{\frac{1}{1+\epsilon}}}+\frac{\sqrt{nb}}{b\sqrt{\log(n)}}\widetilde{\rho}^{\frac{\sqrt{nb}}{\log(n)}}\right).$$ 
Then the result follows from the bandwidth conditions.$\square$

Now, we consider some estimators of $\pi_u$ and $Q_u$. Let 
$\hat{\pi}_u(x)=\sum_{i=1}^{n-1} e_i(u)\mathds{1}_{X_{n,i}=x}$ and $\hat{Q}_u(x,y)=\frac{\hat{\pi}_{u,2}(x,y)}{\hat{\pi}_u(x)}.$
where $\hat{\pi}_{u,2}(x,y)=\sum_{i=1}^{n-1}e_i(u) \mathds{1}_{X_{n,i}=x, X_{n,i+1}=y}$.
\begin{theo}\label{stat}
Assume that for a given $\epsilon>0$, $b\rightarrow 0$ and $nb^{1+\epsilon}\rightarrow \infty$. 
\begin{enumerate}
\item
For $(x,y)\in E^2$, we have
\begin{equation}\label{biais}
\sup_{u\in [0,1]}\left\vert \E\hat{\pi}_u(x)-\pi_u(x)\right\vert=O\left(b^{\kappa}\right),\quad \sup_{u\in [0,1]}\left\vert \frac{\E\hat{\pi}_{u,2}(x,y)}{\E\hat{\pi}_u(x)}-Q_u(x,y)\right\vert=O\left(b^{\kappa}\right)
\end{equation}
and 
For $(x,y)\in E^2$, we have
\begin{equation}\label{variance}
\sup_{u\in [0,1]}\left\vert \hat{\pi}_u(x)-\E\hat{\pi}_u(x)\right\vert=O\left(\frac{\sqrt{\log(n)}}{\sqrt{nb}}\right),\quad \sup_{u\in [0,1]}\left\vert \hat{Q}_u(x,y)-\frac{\E\hat{\pi}_{u,2}(x,y)}{\E\hat{\pi}_u(x)}\right\vert=O\left(\frac{\sqrt{\log(n)}}{\sqrt{nb}}\right).
\end{equation}
\item
For $(u,x)\in [0,1]\times E$, the vector $\left(\sqrt{nb}\left[\hat{\pi}_u(x)-\E\hat{\pi}_u(x)\right]\right)_{x\in E}$ is asymptotically Gaussian with mean $0$ and covariance $\Sigma^{(1)}_u:E\times E\rightarrow \R$ defined by 
$$\Sigma_u^{(1)}=\int K^2(x)dx\cdot \left[\Gamma_u(0)+\sum_{j\geq 1}\left(\Gamma_u(j)+\Gamma_u(j)'\right)\right],$$
where $\Gamma_u(j)_{x,y}=\pi_u(x)Q_u^j(x,y)-\pi_u(x)\pi_u(y)$. 
\item
For $(u,x,y)\in [0,1]\times E^2$,
the vector   
$$\sqrt{nb}\left(\hat{Q}_u(x,y)-\frac{\E\hat{\pi}_{u,2}(x,y)}{\E\hat{\pi}_u(x)}\right)_{(x,y)\in E^2}$$
is asymptotically Gaussian with mean $0$ and covariance $\Sigma^{(2)}:E^2\times E^2\rightarrow \R$ defined by 
$$\Sigma^{(2)}_u\left((x,y),(x',y')\right)=\int K^2(x)dx\cdot \frac{1}{\pi_u(x)}Q_u(x,y)\left[\mathds{1}_{y=y'}-Q_u(x',y')\right]\mathds{1}_{x=x'}.$$  
\end{enumerate}
\end{theo}
\paragraph{Note.} Our estimators are localized versions of the standard estimators used in the homogeneous case. One can see that their convergence rates are standard for nonparametric kernel estimation. 
\paragraph{Proof of Theorem \ref{stat}}
\begin{enumerate}
\item
For the control of the bias, note that 
$$\E\hat{\pi}_{u,2}(x,y)-\pi_{u,2}(x,y)=\sum_{i=1}^{n-1}e_i(u)\left[\pi^{(n)}_{i,2}(x,y)-\pi_{u,2}(x,y)\right].$$
Since $e_i(u)=0$ if $\vert u-i/n\vert>b$, Theorem \ref{first} ensures that 
$$\sup_{u\in [0,1]}\left\vert \E\hat{\pi}_{u,2}(x,y)-\pi_{u,2}(x,y)\right\vert =O\left(b^{\kappa}+\frac{1}{n^{\kappa}}\right)=O\left(b^{\kappa}\right).$$
By summation on $y$, we deduce the first bound in (\ref{biais}) and using the fact that $\min_{u\in [0,1]}\pi_u(x)>0$, we deduce that
$\max_{u\in[0,1]}\frac{1}{\E\hat{\pi}_u(x)}=O_{\P}(1)$ and the second bound in (\ref{biais}) follows.\\
For the variance terms in (\ref{variance}), we use Proposition \ref{exp} which ensures the first bound as well as 
$\max_{u\in[0,1]}\frac{1}{\hat{\pi}_u(x)}=O_{\P}(1)$. This gives also the second bound.
\item
The proof is based on a central limit theorem for triangular arrays of strongly mixing random variables proved in \citet{Rio2}.
This result is given in Proposition \ref{TLC}. For simplicity of notations, we consider the quantity $\sum_{i=1}^ne_i(u)\mathds{1}_{X_{n,i}=x}$ instead of $\hat{\pi}_u(x)$ which has the same asymptotic behavior.
For $x\in E$, let $\lambda_x$ be a real number. We consider the random variables $Z^{(n)}_i=\sum_{x\in E}\lambda_x \mathds{1}_{X_{n,i}=x}$ and $Z_i(u)=\sum_{x\in E} \lambda_x \mathds{1}_{X_i(u)=x}$ and set  
$$G_i^{(n)}=\sqrt{nb}e_i(u)\left(Z^{(n)}_i-\E Z^{(n)}_i\right),\quad H_i^{(n)}=G^{(n)}_i/\sqrt{\v\left(\sum_{j=1}^n G^{(n)}_j\right)}.$$
Let us first derive the limit of $\v\left(\sum_{j=1}^n G^{(n)}_j\right)$.
Using Proposition \ref{mixing1}, we know that there exists a constant $D>0$ and $\widetilde{\rho}\in (0,1)$, such that 
\begin{equation}\label{supercov}
\left\vert\c\left(Z^{(n)}_i,Z^{(n)}_j\right)\right\vert\leq D \widetilde{\rho}^{\vert i-j\vert}.
\end{equation}
Moreover the same type of inequality holds for $\c\left(Z_i(u),Z_j(u)\right)$. 
Then if $\ell$ is a positive integer, let $V_n(\ell)=\left\{(i,j)\in\{1,2,\ldots,n\}^2: \vert i-j\vert\leq \ell\right\}$.  We have 
\begin{eqnarray*}
\v\left(\sum_{j=1}^n G^{(n)}_j\right)&=& \sum_{(i,j)\in V_n(\ell)}\c\left(G^{(n)}_i,G^{(n)}_j\right)+\sum_{(i,j)\in \{1,\ldots,n\}^2\setminus V_n(\ell)} \c\left(G^{(n)}_i,G^{(n)}_j\right)\\
&=& A_1+A_2.
\end{eqnarray*}
If $G_i(u)=\sqrt{nb}e_i(u)\left(Z_i(u)-\E Z_i(u)\right)$, we can also decompose 
\begin{eqnarray*}
\v\left(\sum_{j=1}^n G_j(u)\right)&=& \sum_{(i,j)\in V_n(\ell)}\c\left(G_i(u),G_j(u)\right)+\sum_{(i,j)\in \{1,\ldots,n\}^2\setminus V_n(\ell)} \c\left(G_i(u),G_j(u)\right)\\
&=& A_1(u)+A_2(u).
\end{eqnarray*}
Using (\ref{supercov}), we have 
\begin{eqnarray*}
\vert A_2\vert &\leq& 2nbD\sum_{i=1}^ne_i(u)\sum_{j=i+\ell+1}^n e_j(u) \widetilde{\rho}^{j-i}\\
&\leq& 2nb D \max_{1\leq j\leq n}e_j(u)\frac{\widetilde{\rho}^{\ell}}{1-\rho}.\\
&=& O\left(\widetilde{\rho}^{\ell}\right).
\end{eqnarray*}
In the same way, $\vert A_2(u)\vert=O\left(\widetilde{\rho}^{\ell}\right)$. 
Moreover, using Theorem \ref{first}, we have 
\begin{eqnarray*}
\left\vert A_1-A_1(u)\right\vert&\leq& 2Cnb\sum_{i=1}^ne_i(u)\sum_{j=i}^{(i+\ell)\wedge n}e_j(u)\left[\sum_{s=i}^j\left\vert u-\frac{s}{\ell}\right\vert^{\kappa}+\frac{1}{n^{\kappa}}\right]\\
&\leq& 4 Cnb\max_{1\leq j\leq n}e_j(u)\left[(\ell+1)^2 b^{\kappa}+\frac{\ell+1}{n^{\kappa}}\right]\\
&=& O\left(\ell^2b^{\kappa}+\frac{\ell}{n^{\kappa}}\right).
\end{eqnarray*}
Then, choosing $\ell=\ell_n$ such that $\ell\rightarrow\infty$, $\ell^2 b^{\kappa}\rightarrow 0$ and $\ell/n^{\kappa}\rightarrow 0$, we deduce that
\begin{equation}\label{equiva} 
\v\left(\sum_{j=1}^n G^{(n)}_j\right)=\v\left(\sum_{j=1}^n G_j(u)\right)+o(1).
\end{equation}
Now, we have 
\begin{eqnarray*}
\v\left(\sum_{j=1}^n G_j(u)\right)&=&nb \sum_{i=1}^n e_i(u)^2\v\left(Z_0(u)\right)\\
&+& 2 nb \sum_{i=1}^n\sum_{j=i+1}^n e_i(u)e_j(u)\c\left(Z_0(u),Z_{j-i}(u)\right)\\
&=& nb \sum_{i=1}^n e_i(u)^2\v\left(Z_0(u)\right)\\
&+& 2nb \sum_{s=1}^{n-1}\left[\sum_{i=1}^{n-s} e_i(u)e_{i+s}(u)\right]\sum_{x,y\in E}\lambda_x\lambda_y \Gamma_u(s)_{x,y}.\\
\end{eqnarray*}
Using the Lebesgue theorem and elementary computations with Riemanian sums involving the kernel, we deduce that 
$$\lim_{n\rightarrow \infty}\v\left(\sum_{j=1}^n G_j(u)\right)=\sum_{x,y\in E}\lambda_x\lambda_y \Sigma^{(1)}_{u,x,y}.$$
Using (\ref{equiva}), we also deduce that 
\begin{equation}\label{limitvar} 
\lim_{n\rightarrow \infty}\v\left(\sum_{j=1}^n G_j^{(n)}\right)=\sum_{x,y\in E}\lambda_x\lambda_y \Sigma^{(1)}_{u,x,y}.
\end{equation}
Next, in order to apply Proposition \ref{TLC}, we first check condition (\ref{condvar}). We have $V_{n,n}=1$ and 
$$V_{n,i}=nb\sum_{s,t=1}^i e_s(u)e_t(u)\c\left(Z^{(n)}_s,Z^{(n)}_t\right)\leq nb \sum_{s,t=1}^ne_s(u)e_t(u) \left\vert \c\left(Z^{(n)}_s,Z^{(n)}_t\right)\right\vert= O(1),$$
using (\ref{supercov}). This entails condition (\ref{condvar}) of Proposition \ref{TLC}. 
Finally, we check condition (\ref{condmix}) of Proposition \ref{TLC}. From (\ref{limitvar}), we have $\max_{1\leq i\leq n} \left\vert H^{(n)}_i\right\vert\leq C\mathds{1}_{u-nb\leq i\leq u+nb}/\sqrt{nb}$ for a non random real number $C$ which does not depend on $n$. Then we have also $Q_{n,i}(x)\leq C\mathds{1}_{u-nb\leq i\leq u+nb}/\sqrt{nb}$.
Moreover, $\alpha_{(n)}(x)$ is bounded by (up to a constant) $-\log(x)+1$. This entails that 
$$V_{n,n}^{-3/2}\sum_{i=1}^n\int_0^1 \alpha_{(n)}^{-1}(x/2)Q_{n,i}^2(x)\inf\left(\alpha_{(n)}^{-1}(x/2)Q_{n,i}(x),\sqrt{V_{n,n}}\right)dx=O\left(\frac{1}{\sqrt{nb}}\right).$$
Then we deduce the result of point $2$ from Proposition \ref{TLC}, (\ref{limitvar}) and the Cramér-Wold device.

\item 
Let
$$Z_n(x,y)=\frac{\sqrt{nb}}{\hat{\pi}_u(x)}\sum_{i=1}^{n-1}D_{n,i}(x,y)$$
where
$$D_{n,i}(x,y)=e_i(u)\left[\mathds{1}_{X_{n,i}=x,X_{n,i+1}=y}-Q_{\frac{i+1}{n}}(x,y)\mathds{1}_{X_{n,i}=x}\right]$$
is a martingale increment bounded by $(nb)^{-1}$ (up to a constant). Using the classical Lindeberg central limit theorem 
for martingales, the sum $\sqrt{nb}\sum_{i=1}^{n-1}\left[D_{n,i}(x,y)\right]_{x,y\in E}$ is asymptotically a Gaussian 
vector with mean $0$ and variance matrix $\Sigma$ defined by 
\begin{eqnarray*}
&&\Sigma\left((x,y),(x',y')\right)\\&=&\lim_{n\rightarrow \infty}nb\sum_{i=1}^{n-1}e_i(u)^2\c\left[\mathds{1}_{X_{n,i}=x,X_{n,i+1}=y}-Q_{\frac{i+1}{n}}(x,y)\mathds{1}_{X_{n,i}=x},\mathds{1}_{X_{n,i}=x',X_{n,i+1}=y'}-Q_{\frac{i+1}{n}}(x',y')\mathds{1}_{X_{n,i}=x'}\right]\\
&=& \lim_{n\rightarrow \infty}nb\sum_{i=1}^{n-1}e_i(u)^2\c\left[\mathds{1}_{X_i(u)=x,X_{i+1}(u)=y}-Q_u(x,y)\mathds{1}_{X_i(u)=x},\mathds{1}_{X_i(u)=x',X_{i+1}(u)=y'}-Q_u(x',y')\mathds{1}_{X_i(u)=x'}\right]\\
&=& \int K^2(z)dz \cdot \P\left(X_1(u)=x,X_2(u)=y\right)\cdot\left[\mathds{1}_{y=y'}-Q_u(x',y')\right]\mathds{1}_{x=x'}.
\end{eqnarray*}
In the previous equalities, we have used Theorem \ref{first}, the continuity properties of the transition matrix and the limits 
$$\lim_{n\rightarrow \infty}\frac{1}{nb}\sum_{i=1}^{n-1}K\left(\frac{u-i/n}{b}\right)=\int K(z)dz=1,\quad \lim_{n\rightarrow \infty}\frac{1}{nb}\sum_{i=1}^{n-1}K^2\left(\frac{u-i/n}{b}\right)=\int K(z)^2dz.$$
We deduce that the vector $\left[Z_n(x,y)\right]_{x,y\in E}$ is asymptotically Gaussian with mean zero and covariance matrix $\Sigma_u^{(2)}$.
\\

Then it remains to show that for each $(x,y)\in E^2$,
\begin{equation}\label{final}
\sqrt{nb}\left[\sum_{i=1}^ne_i(u) \frac{\mathds{1}_{X_{n,i}=x}Q_{\frac{i+1}{n}}(x,y)}{\hat{\pi}_u(x)}-\frac{\E\hat{\pi}_{u,2}(x,y)}{\E\hat{\pi}_u(x)}\right]=o_{\P}(1). 
 \end{equation}
To show (\ref{final}), we use the decomposition
\begin{eqnarray*}
&&\sqrt{nb}\left[\sum_{i=1}^{n-1}e_i(u) \frac{\mathds{1}_{X_{n,i}=x}Q_{\frac{i+1}{n}}(x,y)}{\hat{\pi}_u(x)}-\frac{\E\hat{\pi}_{u,2}(x,y)}{\E\hat{\pi}_u(x)}\right]\\
&=&\frac{\sqrt{nb}}{\hat{\pi}_u(x)}\sum_{i=1}^{n-1}e_i(u)\left(\mathds{1}_{X_{n,i}=x}-\pi_i^{(n)}(x)\right)\cdot\left(Q_{\frac{i+1}{n}}(x,y)-Q_u(x,y)\right)\\
&+& \sqrt{nb}\sum_{i=1}^{n-1}e_i(u)\pi_i^{(n)}(x)\left(Q_{\frac{i+1}{n}}(x,y)-Q_u(x,y)\right)\frac{\pi_u(x)-\hat{\pi}_u(x)}{\hat{\pi}_u(x)\pi_u(x)}\\
&=& \frac{A_n}{\hat{\pi}_u(x)}+B_n\frac{\pi_u(x)-\hat{\pi}_u(x)}{\hat{\pi}_u(x)\pi_u(x)}.
\end{eqnarray*} 
Since the kernel $K$ has a compact support and $u\mapsto Q_u(x,y)$ is $\kappa-$Hölder continuous, we have $B_n=O\left(\sqrt{nb}b^{\kappa}\right)$. Moreover, using covariance inequalities, we have $\v(A_n)=O\left(b^{2\kappa}\right)$. 
Then (\ref{final}) follows from $\hat{\pi}_u(x)-\pi_u(x)=O_{\P}\left(\frac{1}{\sqrt{nb}}\right)$ and $\frac{1}{\hat{\pi}_u(x)}=O_{\P}(1)$. The proof of point $3$ is now complete. $\square$ 
\end{enumerate}

\section{Contraction of Markov kernels using Wasserstein metrics}\label{Dob2}
In this section, we consider a Polish space $\left(E,d\right)$. 
For $p\geq 1$,  we consider the set of probability measures on $(E,d)$ admitting a moment of order $p$:
$$\mathcal{P}_p(E)=\left\{\mu\in\mathcal{P}(E): \int d(x,x_0)^p\mu(dx)<\infty\right\}.$$
Here $x_0$ is an arbitrary point in $E$. It is easily seen that the set $\mathcal{P}_p(E)$ does not depend on $x_0$. 

The Wasserstein metric $W_p$ of order $p$ associated to $d$ is defined by 
$$W_p(\mu,\nu)=\inf_{\gamma\in\Gamma(\mu,\nu)}\left\{\int_{E\times E} d(x,y)^pd\gamma(x,y)\right\}^{1/p}$$
where $\Gamma(\mu,\nu)$ denotes the collection of all probability measures on $E\times E$
with marginals $\mu$ and $\nu$. We will say that $\gamma\in \Gamma(\mu,\nu)$ is an optimal coupling of $(\mu,\nu)$ if 
$$\left(\int d(x,y)^p\gamma(dx,dy)\right)^{1/p}=W_p(\mu,\nu).$$
It is well-known that an optimal coupling always exist. See \citet{Villani} for some properties of Wasserstein metrics.  

In the sequel, we will use the following assumptions. 
\begin{description}
\item [B1]
For all $(u,x)\in [0,1]\times E$, $\delta_x Q_u\in\mathcal{P}_p(E)$.

\item[B2]
There exist a positive integer $m$ and two real numbers $r\in (0,1)$ and $C_1\geq 1$ such that for all $u\in [0,1]$ and all $x\in E$,
$$W_p\left(\delta_x Q^m_u,\delta_y Q_u^m\right)\leq rd(x,y),\quad W_p\left(\delta_x Q_u,\delta_y Q_u\right)\leq C_1 d(x,y).$$

\item [B3]
The family of transitions $\left\{Q_u:u\in [0,1]\right\}$ satisfies the following Hölder type continuity condition.
There exist $\kappa\in (0,1]$ and $C_2>0$, such that for all $x\in E$ and all $u,v\in [0,1]$, 
$$W_p\left(\delta_x Q_u,\delta_x Q_v\right)\leq C_2\left(1+d(x,x_0)\right)\vert u-v\vert^{\kappa}.$$
\end{description}

\paragraph{Note.} If $R$ is a Markov kernel, the Dobrushin contraction coefficient is now defined by 
$$c(R)=\sup_{(x,y)\in E^2 \atop x\neq y}\frac{W_p\left(\delta_x R,\delta_y R\right)}{d(x,y)}.$$
Thus Assumption ${\bf B2}$ means that $\sup_{u\in [0,1]}c\left(Q_u\right)<\infty$ and $\sup_{u\in [0,1]}c\left(Q^m_u\right)<1$.

The following proposition shows that under these assumptions, the marginal distribution of the Markov chain with transition $Q_u$ converges exponentially fast to its unique invariant probability distribution which is in turn Hölder continuous with respect to $u$, in Wasserstein metric. 
\begin{prop}\label{inherit}
Assume that assumptions {\bf B1-B3} hold true and set for an integer $j\geq 1$,
\begin{enumerate}
\item
 For all $u\in [0,1]$, the Markov chain of transition $Q_u$ has a unique invariant probability distribution denoted by $\pi_u$. 
Moreover for all initial probability distribution $\mu\in\mathcal{P}_p(E)$, we have for $n=mj+s$
$$W_p\left(\mu Q_u^n,\pi_u\right)\leq C_1^sr^j \left[\left(\int d(x,x_0)^p\mu(dx)\right)^{1/p}+
\kappa_2\right],$$
where $\kappa_2=\sup_{u\in [0,1]}\left(\int d(x,x_0)^p\pi_u(dx)\right)^{1/p}$.
\item
If $u,v\in [0,1]$, we have
$$W_p(\pi_u,\pi_v)\leq \frac{C_2 \vert u-v\vert^{\kappa}}{1-r}\left[m C_1^{m-1}\kappa_2+\sum_{j=0}^{m-1}C_1^j \kappa_1(m-j-1)\right],$$
where $\kappa_1(j)=\sup_{u\in[0,1]}\left(\int d(x,x_0)^p Q^j_u(x_0,dx)\right)^{1/p}$. 
\end{enumerate}
\end{prop}

\paragraph{Proof of Proposition \ref{inherit}}
We first show that the quantities $\kappa_1(j)$ are finite.
We set $q_j=\left(\int \left(1+d(x,x_0)\right)^p Q^j_0(x_0,dx)\right)^{1/p}$.
If $j\geq 1$, we have, using Lemma \ref{preli},
\begin{eqnarray*}
&&W_p\left(\delta_{x_0}Q^j_u,\delta_{x_0}Q^j_0\right)\\
&\leq& W_p\left(\delta_{x_0}Q^j_u,\delta_{x_0}Q^{j-1}_0Q_u\right)+W_p\left(\delta_{x_0}Q^{j-1}_0Q_u,\delta_{x_0}Q^j_0\right)\\
&=& C_1 W_p\left(\delta_{x_0}Q^{j-1}_u,\delta_{x_0}Q^{j-1}_0\right)+C_2\vert u\vert^{\kappa}q_{j-1}. 
\end{eqnarray*}
Then we obtain 
\begin{equation}\label{utility}
W_p\left(\delta_{x_0}Q^j_u,\delta_{x_0}Q^j_0\right)\leq C_2\sum_{s=0}^{j-1}C_1^s q_{j-s-1}.
\end{equation}
Then, using Lemma \ref{Lip} for the function $f(x)=1+d(x,x_0)$, we get 
$$\kappa_1(j)\leq q_j+C_2\sum_{s=0}^{j-1}C_1^s q_{j-s-1}.$$
\begin{enumerate}
\item
The existence and unicity of an invariant probability $\pi_u\in\mathcal{P}_p$ easily follows from the fixed point theorem for a contractant application in the complete metric space $\left(\mathcal{P}_p,W_p\right)$.\\
Before proving the geometric convergence, let us show that the quantity $\kappa_2$ is finite.  
We have, using Lemma \ref{preli},
\begin{eqnarray*}
W_p(\pi_u,\pi_0)&\leq& W_p\left(\pi_uQ^m_u,\pi_0Q^m_u\right)+W_p\left(\pi_0Q^m_u,\pi_0Q^m_0\right)\\
&\leq& r W_p\left(\pi_u,\pi_0\right)+\left(\int W^p_p\left(\delta_x Q^m_u,\delta_x Q^m_0\right)\pi_0(dx)\right)^{1/p}.\\
\end{eqnarray*}
Using (\ref{utility}) and Lemma (\ref{preli}), we have 
\begin{eqnarray*}
W_p\left(\delta_x Q^m_u,\delta_x Q^m_0\right)&\leq & W_p\left(\delta_x Q^m_u,\delta_{x_0} Q^m_u\right)+ W_p\left(\delta_{x_0} Q^m_u,\delta_{x_0} Q^m_0\right)+W_p\left(\delta_x Q^m_0,\delta_{x_0} Q^m_0\right)\\
&\leq& 2 r d(x,x_0)+C_2\sum_{s=0}^{m-1}C_1^j q_{m-s-1}.
\end{eqnarray*}
From the previous bound, we easily deduce the existence of a real number $D>0$, not depending on $u$, such that $W_p(\pi_u,\pi_0)\leq \frac{D}{1-r}$. Then, using Lemma \ref{Lip}, we get 
$$\kappa_2\leq \frac{D}{1-r}+\left(\int d(x,x_0)^p \pi_0(dx)\right)^{1/p},$$
which is finite.\\  
Now, the geometric convergence is a consequence of the inequality
$$W_p\left(\mu Q_u^n,\pi_uQ_u^n\right)\leq C_1^s r^j W_p(\mu,\pi_u)\leq C_1^s r^j\left[\left(\int d(x,x_0)^p\mu(dx)\right)^{1/p}+
\kappa_2\right].$$
Finally, let $\nu$ be an invariant probability for $P_u$ (not necessarily in $\mathcal{P}_p$).
Let $f:E\rightarrow \R$ be an element of $\mathcal{C}_b(E)$. Since convergence in Wasserstein metric implies weak convergence, we have from the geometric ergodicity $\lim_{n\rightarrow \infty}Q^n_u f(x)=\pi_u f$ for all $x\in E$. Hence, using the Lebesgue theorem, we have
$$\nu f =\nu Q_u^n f=\int \nu(dx)Q_u^nf(x)\rightarrow \pi_u f$$
which shows the unicity of the invariant measure. 

\item
Proceeding as for the previous point, we have 
\begin{equation}\label{utility2}
W_p(\pi_u,\pi_v)\leq r W_p(\pi_u,\pi_v)+\left(\int W_p^p\left(\delta_x Q_u^m,\delta_x Q_v^m\right)\pi_v(dx)\right)^{1/p}.
\end{equation}
But 
\begin{eqnarray*}
W_p\left(\delta_x Q_u^m,\delta_x Q_v^m\right)&\leq& C_1 W_p\left(\delta_x Q_u^{m-1},\delta_x Q_v^{m-1}\right)+C_2\vert u-v\vert^{\kappa}\left(\int \left[1+d(y,x_0)\right]^p Q^{m-1}_v(x,dy)\right)^{1/p}\\
&\leq &C_1 W_p\left(\delta_x Q_u^{m-1},\delta_x Q_v^{m-1}\right)+C_2\vert u-v\vert^{\kappa}\left(\kappa_1(m-1)+C_1^{m-1} d(x,x_0)\right).
\end{eqnarray*}
We deduce that 
$$W_p\left(\delta_x Q_u^m,\delta_x Q_v^m\right)\leq C_2\vert u-v\vert^{\kappa}\left(\sum_{j=0}^{m-1}C_1^j\kappa_1(m-j-1)+mC_1^{m-1} d(x,x_0)\right).$$
Reporting the last bound in (\ref{utility2}), we get the result.$\square$ 
\end{enumerate}

Now let us give the main result of this section. For $j\in\N^{*}$, we endow the space $E^j$ with the distance
$$d_j(x,y)=\left(\sum_{s=1}^j d(x_s,y_s)^p\right)^{1/p},\quad x,y\in E^j.$$
We will still denote by $W_p$ the Wasserstein metric for Borelian measures on $E^j$. 

\begin{theo}\label{central}
Assume that assumptions ${\bf B1-B3}$ hold true. 
Then the triangular array of Markov chains $\left\{X_{n,k}: n\in\Z^+, k\leq n\right\}$ is locally stationary.
Moreover, there exists a real number $C>0$,
only depending on $j,p,d,\kappa,r,C_1,C_2,\kappa_1(1),\ldots,\kappa_1(m),\kappa_2$ such that
$$W_p\left(\pi^{(n)}_{k,j},\pi_{u,j}\right)\leq C\left[\sum_{s=k}^{k+j-1}\left\vert u-\frac{s}{n}\right\vert^{\kappa}+\frac{1}{n^{\kappa}}\right].$$
\end{theo}

\paragraph{Proof of Theorem \ref{central}}
\begin{enumerate}
\item
We show the result by induction and first consider the case $j=1$. For $k\leq n$, 
let $Q_{k,m}$ be the probability kernel $Q_{\frac{k-m+1}{n}}\cdots Q_{\frac{k}{n}}$. 
We have 
\begin{eqnarray*}
W_p\left(\pi_k^{(n)},\pi_u\right)&=&W_p\left(\pi_{k-m}^{(n)}Q_{k,m},\pi_u Q_u^m\right)\\
&\leq& W_p\left(\pi_{k-m}^{(n)}Q_{k,m},\pi_{k-m}^{(n)}Q_u^m\right)+W_p\left(\pi_{k-m}^{(n)}Q_u^m,\pi_u Q_u^m\right)\\
&\leq & r W_p\left(\pi_{k-m}^{(n)},\pi_u\right)+\left(\int W_p^p\left(\delta_x Q_{k,m},\delta_x Q_u^m\right)\pi_{k-m}^{(n)}(dx)\right)^{1/p}.
\end{eqnarray*}
From Lemma \ref{preli2}, we have 
$$W_p\left(\delta_x Q_{k,m},\delta_x Q_u^m\right)\leq \sum_{s=0}^{m-1}C_1^{s}C_2\left[\int \left(1+d(y,x_0)\right)^p \delta_x Q_{\frac{k-m+1}{n}}\cdots Q_{\frac{k-s-1}{n}}(dy)\right]^{1/p}\left\vert u-\frac{k-s}{n}\right\vert^{\kappa}.$$
First we note that from our assumptions and using Lemma \ref{Lip} for the function $f(x)=1+d(x,x_0)$, we have 
$$\left[\delta_xQ_u f^p\right]^{1/p}\leq \left[\delta_{x_0}Q_uf^p\right]^{1/p}+C_1 d(x,x_0)\leq \left(1+\kappa_1(1)+C_1\right)f(x),$$
where $\kappa_1$ is defined in Proposition \ref{inherit}. Then we get $\sup_{u\in [0,1]}\delta_xQ_u f^p\leq C_3^p f^p(x)$,
where $C_3=1+\kappa_1(1)+C_1$. 
This yields to the inequality 
$$W_p\left(\delta_x Q_{k,m},\delta_x Q_u^m\right)\leq \sum_{s=0}^{m-1}C_1^{s}C_2 C_3^{m-s-1}\left\vert u-\frac{k-s}{n}\right\vert^{\kappa}f(x).$$
Then we obtain 
$$W_p\left(\pi_k^{(n)},\pi_u\right)\leq  r W_p\left(\pi_{k-m}^{(n)},\pi_u\right)+\sum_{s=0}^{m-1}C_1^{s}C_2C_3^{m-s-1}\left\vert u-\frac{k-s}{n}\right\vert^{\kappa}\left(\pi_{k-m}^{(n)}f^p\right)^{1/p}.$$
Then the result will easily follow if we prove that $\sup_{n,k\leq n}\pi^{(n)}_k f^p<\infty$.
Setting $c_k=W_p\left(\pi_k^{(n)},\pi_{\frac{k}{n}}\right)$ and $C_4=\sum_{s=0}^{m-1}(s+1)^{\kappa}C_1^{s+1}C_2C_3^{m-s-1}$ and using our previous inequality, we have 
\begin{eqnarray*}
c_k&\leq & r W_p\left(\pi_{k-m}^{(n)},\pi_{\frac{k}{n}}\right)+\frac{C_4}{n^{\kappa}}\left(\pi_{k-m}^{(n)}f^p\right)^{1/p}\\
&\leq& \left(r+\frac{C_4}{n^{\kappa}}\right)c_{k-m}+rW_p\left(\pi_{\frac{k-m}{n}},\pi_{\frac{k}{n}}\right)+\frac{C_4}{n^{\kappa}}(1+\kappa_2).
\end{eqnarray*}
Then, if $n_0$ is such that for all $n\geq n_0$, $r+\frac{C_4}{n^{\kappa}}<1$, the last inequality, Proposition \ref{inherit} and Lemma \ref{Lip} guarantee that 
$\sup_{n\geq n_0,k\leq n}\pi^{(n)}_k f^p$ is finite and only depends on $p,d,r,C_1,C_2,\kappa_1(1),\ldots,\kappa_1(m),\kappa_2,\kappa$.
Moreover if $n\leq n_0$, we have  $\left(\pi^{(n)}_k f^p\right)^{1/p}\leq (C_4+1)^{n_0}\left(\pi_0 f^p\right)^{1/p}$. This concludes the proof for the case $j=1$. 
\item
Now for $j\geq 2$, we define a coupling of $\left(\pi^{(n)}_{k,j},\pi_{u,j}\right)$ as follows. First we consider an optimal coupling $\Gamma_{u,j-1}^{(k,n)}$ of $\left(\pi^{(n)}_{k,j-1},\pi_{u,j-1}\right)$, and for each $(x,y)\in E^2$, we define an optimal coupling 
$\Delta_{x,y,j,u}^{(k,n)}$ of $\left(\delta_x Q_{\frac{k+j}{n}},\delta_y Q_u\right)$. From \citet{Villani}, Corollary $5.22$, it is possible to choose this optimal coupling such that the application $(x,y)\mapsto \Delta_{x,y,j,u}^{(k,n)}$ is measurable. Now we define
$$\Gamma_{u,j}^{(k,n)}(dx_1,dy_1,\ldots,dx_j,dy_j)=\Delta_{x_{j-1},y_{j-1},j,u}^{(k,n)}(dx_j,dy_j) \Gamma_{u,j-1}^{(k,n)}(dx_1,dy_1,\ldots,dx_{j-1},dy_{j-1}).$$
Then we easily deduce that 
$$W_p^p\left(\pi^{(n)}_{k,j},\pi_{u,j}\right)\leq W_p^p\left(\pi^{(n)}_{k,j-1},\pi_{u,j-1}\right)+\int W_p^p\left(\delta_{x_{j-1}}Q_{\frac{k+j}{n}},\delta_{y_{j-1}}Q_u\right)\Gamma^{(n,k)}_{u,j-1}\left(dx_1,dy_1,\ldots,dx_{j-1},dy_{j-1}\right).$$
Since 
$$W_p\left(\delta_{x_{j-1}}Q_{\frac{k+j}{n}},\delta_{y_{j-1}}Q_u\right)\leq C_1 d(x_{j-1},y_{j-1})+C_2\left[1+d(y_{j-1},x_0)\right]\left\vert u-\frac{k+j}{n}\right\vert^{\kappa}.$$
This leads to 
$$W_p\left(\pi^{(n)}_{k,j},\pi_{u,j}\right)\leq (1+C_1)W_p\left(\pi^{(n)}_{k,j-1},\pi_{u,j-1}\right)+C_2(1+\kappa_2)\left\vert u-\frac{k+j}{n}\right\vert^{\kappa}.$$
The results follows by a finite induction.

Finally, note that Condition $1$ of Definition \ref{debase} follows from induction and the point $2$ of Proposition \ref{inherit}, because using the same type of arguments, we have 
$$W_p\left(\pi_{v,j},\pi_{u,j}\right)\leq (1+C_1)W_p\left(\pi_{v,j-1},\pi_{u,j-1}\right)+C_2(1+\kappa_2)\left\vert u-v\right\vert^{\kappa}.$$

The proof of the Theorem is now complete.$\square$
\end{enumerate}

\subsection{Mixing conditions}
We now introduce another useful coefficient: the $\tau-$mixing coefficient introduced and studied in \citet{ded} that we will adapt to our triangular arrays. This coefficient
has been introduced for Banach spaces $E$. In the sequel, we denote by $\Lambda_1(E)$ the set of $1-$Lipschitz functions from $E$ to $\R$. Assume first that $E=\R$ and as for the $\beta-$mixing coefficients, set $X_{n,j}=0$ for $j>n$. 
Then setting 
$$U\left(\mathcal{F}^{(n)}_i,X_{n,i+j}\right)=\sup\left\{\left\vert \E\left[f\left(X_{n,i+j}\vert \mathcal{F}^{(n)}_i\right)\right]-\E\left[f\left(X_{n,i+j}\right)\right]\right\vert:\quad f\in \Lambda_1(\R)\right\},$$
the $\tau_n-$mixing coefficient for the sequence $\left(X^{(n)}_j\right)_{j\in\Z}$  
is defined by 
$$\tau_n(j)=\sup_{i\in\Z}\E\left[U\left(\mathcal{F}^{(n)}_i,X_{n,i+j}\right)\right].$$
Now for a general metric space $E$, the $\tau_n-$mixing coefficient is defined by  
$$\tau_n(j)=\sup_{f\in \Delta_1(E)}\sup_{i\in\Z}\E\left[U\left(\mathcal{F}^{(n)}_i,f\left(X_{n,i+j}\right)\right)\right].$$
Note that, if $\widetilde{X}_{n,i}$ denotes a copy of $X_{n,i}$,
$$\tau_n(j)\leq \sup_{i\in\Z}\E d\left(X_{n,i},\widetilde{X}_{n,i}\right),$$

For bounding this mixing coefficient, the following assumption, which strengthens assumption {\bf B2} in the case $m\geq 2$ and $p=1$, will be needed.
\begin{description}
\item
{\bf B4} There exists a positive real number $\epsilon$ such that for all $\left(u,u_1,\ldots,u_m\right)\in [0,1]^{m+1}$ satisfying $\vert u_i-u\vert<\epsilon$ for $1\leq i\leq m$, we have
$$W_1\left(\delta_xQ_{u_1}\cdots Q_{u_m},\delta_yQ_{u_1}\cdots Q_{u_m}\right)\leq r d(x,y),$$
where $m$ and $r$ are defined in assumption {\bf B2}.
\end{description}

\begin{prop}\label{mixing2}
Assume that assumptions ${\bf B2}$ and ${\bf B4}$ hold true. 
Then there exists $C>0$ and $\rho\in (0,1)$, only depending on $m,r,C_1,\epsilon$ such that 
$$\tau_n(j)\leq C \rho^j.$$
\end{prop}

\paragraph{Proof of Proposition \ref{mixing2}}
We first consider the case $n\geq m/\epsilon$. Now if $k$ is an integer such that $k+m-1\leq n$, note that assumption {\bf B4} entails that 
\begin{equation}\label{mimix2}
W_1\left(\mu Q_{\frac{k}{n}}\cdots Q_{\frac{k+m-1}{n}},\nu Q_{\frac{k}{n}}\cdots Q_{\frac{k+m-1}{n}}\right)\leq r W_1(\mu,\nu),
\end{equation}
where the probability measures $\mu$ and $\nu$ have both a finite first moment. 
If $j=m t+s$, we get from (\ref{mimix2}) and Assumption ${\bf B2}$,
$$\tau_n(j)\leq C_1^sr^t \sup_{i\in\Z}\E\left[W_1\left(\delta_{X_{n,i}},\pi^{(n)}_i\right)\right]\leq  2 \sup_{i\in\Z}\E d\left( X_{n,i},x_0\right)\cdot C_1^s r^t.$$
We have seen in the proof of Theorem \ref{central} that $\sup_{n\in\Z,i\leq n}\E d\left( X_{n,i},x_0\right)<\infty$. 

Now assume that $n<m/\epsilon$. If $j\leq n$, we have 
$$\tau_n(j)\leq 2 \sup_{i\in\Z}\E d\left( X_{n,i},x_0\right)\cdot C_1^{m/\epsilon}.$$
Now if $j>n$, we have since $(X_{n,j})_{j\leq 0}$ is stationary with transition kernel $Q_0$, 
$$\tau_n(j)\leq 2 \sup_{i\in\Z}\E d\left( X_{n,i},x_0\right)\cdot C_1^{m/\epsilon+m}r^{\left[\frac{j-n}{m}\right]}.$$
This leads to the result for $\rho=r^{1/m}$ and an appropriate choice of $C>0$.$\square$

\paragraph{Note.} Let us remind that Proposition \ref{mixing2} implies a geometric decrease for the covariances. This is a consequence of the following property. If 
$f:E\rightarrow \R$ is measurable and bounded and $g:E\rightarrow \R$ is measurable and Lipschitz, we have 
$$\c\left(f\left(X_{n,i}\right),g\left(X_{n,i+j}\right)\right)\leq \Vert f\Vert_{\infty}\cdot\delta(g)\cdot\tau_n(j).$$

\subsection{An extension to $q-$order Markov chains}

We start with an extension of our result to Markov sequences of order $q\geq 1$ and taking values in the Polish space $(E,d)$.
Let $\left\{S_u: u\in [0,1]\right\}$ be a family of probability kernels from $\left(E^q,\mathcal{B}(E^q)\right)$ to $\left(E,\mathcal{B}(E)\right)$. The two following assumptions will be used.

\begin{description}
\item [H1]
For all ${\bf x}\in E^q$, $S_u({\bf x},\cdot)\in\mathcal{P}_p(E)$.
\item [H2]
There exist non-negative real numbers $a_1,a_2,\ldots,a_q$ satisfying $\sum_{j=1}^q a_j<1$ and such that for all $(u,{\bf x},{\bf y})\in [0,1]\times E^q\times E^q$, 
$$W_p\left(S_u({\bf x},\cdot),S_u ({\bf y},\cdot)\right)\leq \sum_{j=1}^q a_j d(x_j,y_j).$$
\item[H2]
There exists a positive real number $C$ and $\kappa\in (0,1)$ such that for all $(u,,v,{\bf x})\in [0,1]\times[0,1]\times E^q$, 
$$W_p\left(S_u({\bf x},\cdot),S_v({\bf x},\cdot)\right)\leq C\left(1+\sum_{j=1}^q d(x_j,x_0)\right)\vert u-v\vert^{\kappa}.$$
\end{description}
To define Markov chains, we consider the family of Markov kernels $\left\{Q_u: u\in [0,1]\right\}$ on the measurable space $\left(E^q,\mathcal{B}(E)^q\right)$ and defined by 
$$Q_u\left({\bf x},d{\bf y}\right)=S_u\left({\bf x},dy_q\right)\otimes \delta_{x_2}(y_1)\otimes\cdots \otimes \delta_{x_q}(d y_{q-1}).$$

\begin{cor}\label{autoreg}
If the assumptions {\bf H1-H3} hold true then Theorem \ref{central} and Proposition \ref{mixing2} apply.
\end{cor}

\paragraph{Proof of Corollary \ref{autoreg}}
Assumption ${\bf H1}$ entails ${\bf B1}$.
Then we check assumption ${\bf B3}$. If $(u,v,{\bf x})\in [0,1]\times [0,1]\times E^q$, let $\alpha_{{\bf x},u,v}$ be a coupling of the two probability distributions $S_u({\mbox x},\cdot)$ and $S_v({\mbox x},\cdot)$. Then 
$$\gamma_{{\bf x},u,v}(d{\bf y},d{\bf y'})=\alpha_{{\bf x},u,v}(dy_q,dy'_q)\otimes_{j=1}^{q-1}\delta_{x_{j+1}}(d y_j)\otimes \delta_{x_{j+1}}(d y'_j)$$
defines a coupling of the two measures $\delta_{{\bf x}}Q_u$ and $\delta_{{\bf x}}Q_v$. We have
$$W_p\left(\delta_{{\bf x}}Q_u,\delta_{{\bf x}}Q_v\right)\leq \left[\int d(y_q,y'_q)^p\alpha_{{\bf x},u,v}(dy_q,dy'_q)\right]^{1/p}.$$  
By taking the infinimum of the last bound over all the couplings, we get 
$$W_p\left(\delta_{{\bf x}}Q_u,\delta_{{\bf x}}Q_v\right)\leq W_p\left(S_u\left({\bf x},\cdot\right),S_v\left({\bf x},\cdot\right)\right),$$
which shows  ${\bf B3}$, using assumption {\bf H3}.\\
Finally, we check assumptions ${\bf B2}$ and ${\bf B4}$.
For an integer $m\geq 1$, $(u_1,\ldots,u_m)\in [0,1]^m$ and $({\bf x},{\bf y})\in E^q\times E^q$, we denote by $\alpha_{{\bf x,y},u}$ an optimal coupling of $\left(S_u({\bf x},\cdot),S_u({\bf y},\cdot)\right)$. From \citet{Villani}, Corollary $5.22$, there exists a measurable choice of $({\bf x},{\bf y})\mapsto \alpha_{{\bf xy},u}$. 
We define 
$$\gamma^{({\bf x_{q+1}y_{q+1}})}_{m,u_1,\ldots,u_m}\left(dx_{q+1},\ldots,dx_{q+m},dy_{q+1},\ldots,dy_{q+m}\right)=\prod_{i=q+1}^{q+m} \alpha_{{\bf x_i,y_i},u_i}(d x_i,dy_i),$$
where ${\bf x_i}=(x_{i-1},\ldots,x_{i-q})$. Let $\Omega=E^m\times E^m$ endowed with its Borel sigma field and the probability measure $\P=\gamma^{({\bf x_{q+1}y_{q+1}})}_{m,u_1,\ldots,u_m}$. Then we define the random variables $Z_j^{{\bf x_{q+1}}}=x_j$, $Z_j^{{\bf y_{q+1}}}=y_j$ for $1\leq j\leq q$ and for $1\leq j\leq m$, $Z^{{\bf x_{q+1}}}_{q+j}(\omega_1,\omega_2)=\omega_{1,j}$,
$Z^{{\bf y_{q+1}}}_{q+j}(\omega_1,\omega_2)=\omega_{2,j}$ for $j=1,\ldots,m$.
By definition of our couplings, we have 
$$\E^{1/p}\left[d\left(Z_k^{{\bf x_{q+1}}},Z_k^{{\bf y_{q+1}}}\right)^p\right]\leq \sum_{j=1}^q a_j \E^{1/p}\left[d\left(Z_{k-j}^{{\bf x_{q+1}}},Z_{k-j}^{{\bf y_{q+1}}}\right)^p\right].$$
Using a finite induction, we obtain 
$$\E^{1/p}\left[d\left(Z_k^{{\bf x_{q+1}}},Z_k^{{\bf y_{q+1}}}\right)^p\right]\leq \alpha^{\frac{k}{q}}\max_{1\leq j\leq q}d(x_j,y_j),$$
where $\alpha=\sum_{j=1}^q a_j$.
Setting $X_m^{{\bf x}}=\left(Z^{{\bf x}}_{m-q+1},\ldots,Z^{{\bf x}}_m\right)$, this entails 
\begin{eqnarray*}
W_p\left(\delta_{{\bf x}}Q_{u_1}\cdots Q_{u_m},\delta_{{\bf y}}Q_{u_1}\cdots Q_{u_m}\right)&\leq& \E^{1/p}\left[d_q\left(X_m^{{\bf x}},X_m^{{\bf y}}\right)^p\right]\\
&\leq & \sum_{j=1}^q \alpha^{\frac{m-j+1}{q}}\max_{1\leq j\leq q}d(x_j,y_j)\\
&\leq & \sum_{j=1}^q \alpha^{\frac{m-j+1}{q}}\cdot d_q({\bf x},{\bf y}). 
\end{eqnarray*}
Then {\bf B2-B4} are satisfied if $m$ is large enough by noticing that $W_1\leq W_p$. $\square$

\subsection{Examples of locally stationary Markov chains}

Natural examples of a $q-$order Markov chain satisfying the previous assumptions are based on time-varying autoregressive process. More precisely, if $E$ and $G$ are measurable spaces and $F:[0,1]\times E^q\times G\rightarrow E$, the triangular array $\left\{X_{n,i}: i\leq n, n\in \Z^+\right\}$ is defined recursively by the equations 
\begin{equation}\label{autoreg}
X_{n,i}=F\left(\frac{i}{n},X_{n,i-1},\ldots,X_{n,i-q},\varepsilon_i\right),\quad i\leq n,
\end{equation}
where the usual convention is to assume that 
$$X_{n,i}=F\left(0,X_{n,i-1},\ldots,X_{n,i-q},\varepsilon_i\right),\quad i\leq 0.$$
Then, if $S_u({\bf x},\cdot)$ denotes the distribution of $F\left(u,x_q,\ldots,x_1,\varepsilon_1\right)$, we have 
$$W_p\left(S_u({\bf x},\cdot),S_u ({\bf y},\cdot)\right)\leq \E^{1/p}\left[d\left(F\left(u,x_q,\ldots,x_1,\varepsilon_1\right),F\left(u,y_q,\ldots,y_1,\varepsilon_1\right)\right)^p\right],$$
$$W_p\left(S_u({\bf x},\cdot),S_v ({\bf x},\cdot)\right)\leq \E^{1/p}\left[d\left(F\left(u,x_q,\ldots,x_1,\varepsilon_1\right),F\left(v,x_q,\ldots,x_1,\varepsilon_1\right)\right)^p\right].$$
Then the assumptions ${\bf H1-H3}$ are satisfied if for all $(u,v,{\bf x},{\bf y})\in [0,1]\times[0,1]\times E^q\times E^q$,
$$\E^{1/p}\left[d\left(F(u,x_q,\ldots,x_1),x_0\right)^p\right]<\infty,$$
$$\E^{1/p}\left[d\left(F(u,x_q,\ldots,x_1),F(u,y_q,\ldots,y_1)\right)^p\right]\leq \sum_{j=1}^q a_j d\left(x_{q-j+1},y_{q-j+1}\right)$$
and 
$$\E^{1/p}\left[d\left(F(u,x_q,\ldots,x_1),F(v,x_q,\ldots,x_1)\right)^p\right]\leq C\left(1+\sum_{j=1}^q d(x_j,x_0)\right)\cdot\vert u-v\vert^{\kappa}.$$ 
A typical example of such time-varying autoregressive process is the univariate tv-ARCH process for which 
$$X_{n,i}=\xi_i\sqrt{a_0(i/n)+\sum_{j=1}^q a_j(i/n)X_{n,i-j}^2},$$
with $\E\xi_t=0$, $\v\xi_t=1$. The previous assumptions are satisfied for the square of this process if the $a_j$'s are $\kappa-$Hölder continuous and if 
$$\Vert\xi^2_t\Vert_p\cdot\sup_{u\in [0,1]}\sum_{j=1}^q a_j(u)<1,\mbox{ for some } p\geq 1.$$
See \citet{Fryz} and \citet{Truq} for the use of those processes for modeling financial data.

\paragraph{Note.\label{gen}} The approximation of time-varying autoregressive processes by stationnary processes is discussed in several papers. See for instance \citet{SR} for linear autoregressions with time varying random coefficients, \citet{Vogt} for nonlinear time-varying autoregressions or \citet{ZW} for additional results in the same setting. 
However, the approximating stationary process of (\ref{autoreg}) is given by 
$$X_i(u)=F\left(u,X_{i-1}(u),\ldots,X_{i-q}(u),\varepsilon_i\right).$$
Note that $W_p\left(\pi^{(n)}_k,\pi_u\right)\leq \E^{1/p}\left[d\left(X_k^{(n)},X_k(u)\right)^p\right]$ and the aforementioned references usually study a control of this upper bound by $\left\vert u-\frac{k}{n}\right\vert^{\kappa}+\frac{1}{n^{\kappa}}$.
Note that in the case of autoregressive processes, a coupling of the time-varying processes and its stationary approximation
is already defined because the same noise process is used in both cases. However it is possible to construct some examples for which $\pi_k^{(n)}=\pi_u$ and $\E^{1/p}\left[d\left(X_{n,k},X_k(u)\right)^p\right]\neq 0$, i.e the coupling used is not optimal. 
Nevertheless, it is still possible to obtain an upper bound of $\E^{1/p}\left[d\left(X_k^{(n)},X_k(u)\right)^p\right]$ using our results.
To this end, let us assume that $q=1$ (otherwise one can use vectors of $q-$successive coordinates to obtain a Markov chain of order $1$) if $u\in[0,1]$, and we consider the Markov kernel form $\left(E^2,\mathcal{B}(E^2)\right)$ to itself, given by 
$$Q_v^{(u)}\left(x_1,x_2,A\right)=\P\left(\left(F(v,x_1,\varepsilon_1),F(u,x_2,\varepsilon_1)\right)\in A\right),\quad A\in\mathcal{B}(E^2),\quad v\in [0,1].$$
One can show that the family $\left\{Q_v^{(u)}: v\in [0,1]\right\}$ satisfies the assumptions ${\bf B1-B4}$ for the metric
$$d_2\left[(x_1,x_2),(y_1,y_2)\right]=\left(d(x_1,y_1)^p+d(x_2,y_2)^p\right)^{1/p}.$$ 
Moreover, the constant in Theorem \ref{central} does not depend on $u\in [0,1]$. 
Then Lemma \ref{astuce} guarantees that there exists a positive constant $C$ not depending on $k,n,u$ such that 
$$\E^{1/p}\left(d\left(X_{k,n},X_k(u)\right)^p\right)\leq C\left[\left\vert u-\frac{k}{n}\right\vert^{\kappa}+\frac{1}{n^{\kappa}}\right].$$

\paragraph{Iteration of random affine functions} 
Here we assume that for each $u\in [0,1]$,  there exists a sequence $\left(A_t(u),B_t(u)\right)_{t\in\Z}$ of i.i.d random variables such that $A_t(u)$ takes its values in the space $\mathcal{M}_d$ of squares matrices of dimension $d$ with real coefficients and $B_t(u)$ takes its values in $E=\R^d$. Let $\Vert\cdot\Vert$ a norm on $E$. We also denote by $\Vert\cdot\Vert$ the corresponding operator norm on $\mathcal{M}_d$. We then consider the following recursive equations 
\begin{equation}\label{saig}  
X_{n,i}=A_i\left(\frac{i}{n}\right)X_{n,i-1}+B_i\left(\frac{i}{n}\right). 
\end{equation}
Local approximation of these autoregressive processes by their stationary versions $X_t(u)=A_t(u)X_{t-1}(u)+B_t(u)$ is studied
is studied by \citet{SR}. In this subsection, we will derive similar results using our Markov chain approach.  
For each $u\in [0,1]$, we denote by $\gamma_u$ the top Lyapunov exponent of the sequence $\left(A_t(u)\right)_{t\in\Z}$, i.e
$$\gamma_u=\inf_{n\geq 1}\frac{1}{n}\E\log\Vert A_n(u)A_{n-1}(u)\cdots A_1(u)\Vert.$$
We assume that there exists $t\in (0,1)$ such that
\begin{description}
\item[R1]  
for all $u\in [0,1]$, $\E\Vert A_1(u)\Vert^t<\infty$, $\E\Vert B_1(u)\Vert^t<\infty$ and $\gamma_u<0$. 
\item[R2]
There exists $C>0$ and $\widetilde{\kappa}\in (0,1)$ such that for all $(u,v)\in [0,1]^2$, 
$$\E\Vert A_1(u)-A_1(v)\Vert^t +\E\Vert B_1(u)-B_1(v)\Vert^t\leq C\vert u-v\vert^{t\widetilde{\kappa}}.$$
\end{description}

\begin{prop}\label{fractional}
For $s\in (0,1)$, we set $d(x,y)=\Vert x-y\Vert^s$. Assume that assumptions ${\bf R1-R2}$ hold true. 
Then there exists $s\in (0,t)$ such that Theorem \ref{central} applies with $p=1$, $\kappa=t\widetilde{\kappa}$ and $d(x,y)=\Vert x-y\Vert^s$ and $x_0=0$.
\end{prop}

\paragraph{Notes}
\begin{enumerate}
\item
 Using the remark in the Note of Section \ref{gen}, we also have $\E \Vert X_{n,k}-X_k(u)\Vert^s\leq C\left(\left\vert
u-\frac{k}{n}\right\vert^{\widetilde{\kappa}}+\frac{1}{n^{\widetilde{\kappa}}}\right)^s$, where the process $\left(X_j(u)\right)_{j\in\Z}$ satisfies the iterations $X_k(u)=A_k(u)X_{k-1}(u)+B_k(u)$. Then the triangular array  $\left\{X_{n,k}: k\leq n,n\in\Z^+\right\}$ is locally stationary in the sense given in \citet{Vogt} (see Definition $2.1$ of that paper). 
\item
One can also give additional results for the Wasserstein metric of order $p\geq 1$ and $d(x,y)=\Vert x-y\Vert$ if 
$$\E\Vert A_1(u)\Vert^p+\E\Vert B_1(u)\Vert^p<\infty,\quad \E^{1/p}\Vert A_1(u)-A_1(v)\Vert^p+\E^{1/p}\Vert B_1(u)-B_1(v)\Vert^p\leq C\vert u-v\vert^{\kappa}$$
and there exists an integer $m\geq 1$ such that $\sup_{u\in [0,1]}\E\Vert A_m(u)\cdots A_1(u)\Vert^p<1$. 
In particular, one can recover results about the local approximation of tv-AR processes defined by 
$$X_{n,i}=\sum_{j=1}^q a_j(i/n)X_{n,i-j}+\sigma(i/n)\varepsilon_i$$
by vectorizing $q$ successive coordinates and assuming $\kappa-$Hölder continuity for the $a_j$'s and $\sigma$.
Details are omitted.

\end{enumerate}

\paragraph{Proof of Proposition \ref{fractional}}
For all $(x,u)\in \R^d\times [0,1]$, the measure $\delta_x Q_u$ is the probability distribution of the random variable $A_k(u)x+B_k(u)$.
Condition ${\bf A1}$ of Theorem \ref{central} follows directly from assumption ${\bf R1}$ (whatever the value of $s\in (0,t)$).
Moreover, we have for $s\in (0,t)$,
\begin{eqnarray*}
W_1\left(\delta_x Q_u,\delta_x Q_v\right)&\leq& \E\Vert A_k(u)-A_k(v)\Vert^s\cdot \Vert x\Vert^s+\E\Vert B_k(u)-B_k(v)\Vert^s\\
&\leq& \left(1+\Vert x\Vert^s\right)\cdot\left(\E^{\frac{s}{t}}\Vert A_k(u)-A_k(v)\Vert^t+\E^{\frac{s}{t}}\Vert B_k(u)-B_k(v)\Vert^t\right).
\end{eqnarray*}
This entails condition ${\bf A3}$, using assumption ${\bf R2}$.
Next, if $u\in[0,1]$, the conditions $\gamma_u<0$ and $\E\Vert A_t(u)\Vert^t<\infty$ entail the existence of an integer $k_u$ and $s_u\in (0,t)$ such that $\E\Vert A_{k_u}(u)A_{k_u-1}(u)\cdots A_1(u)\Vert^{s_u}<1$ (see for instance \citet{FZ}, Lemma $2.3$).
Using the axiom of choice, let us select for each $u$, a couple $(k_u,s_u)$ satisfying the previous property.
From assumption ${\bf R2}$, the set 
$$\mathcal{O}_u=\left\{v\in [0,1]: \E\Vert A_{k_u}(v)A_{k_u-1}(v)\cdots A_1(v)\Vert^{s_u}<1\right\}$$
is an open set of $[0,1]$. By a compactness argument, there exist $u_1,\ldots,u_d\in [0,1]$ such that 
$[0,1]=\cup_{i=1}^d\mathcal{O}_{u_i}$. Then setting $s=\min_{1\leq i\leq d}s_{u_i}$ and denoting by $m$ 
the lowest common multiple of the integers $k_{u_1},\ldots,k_{u_d}$, we have from assumption ${\bf R2}$,
$$r=\sup_{u\in [0,1]}\E \Vert A_m(u)\cdots A_1(u)\Vert^s<1.$$
This entails condition ${\bf B2}$ for this choice of $s$, $m$ and $r$. Indeed, we have 
$$W_1\left(\delta_x Q_u^m,\delta_y Q_u^m\right)\leq \E\Vert A_m(u)\cdots A_1(u)(x-y)\Vert^s\leq r d(x,y).$$
Note also that condition ${\bf B4}$ easily follows from the uniform continuity of the application $(u_1,\ldots,u_m)\mapsto \E \Vert A_m(u_1)\cdots A_1(u_m)\Vert^s$
.$\square$.

\paragraph{Time-varying integer-valued autoregressive processes (tv-INAR)}
Stationary INAR processes are widely used in the time series community. This time series model has been proposed by \citet{Al} and a generalization to several lags was studied in \citet{Guan}. 
In this paper, we introduce a locally stationary version of such processes. 
For $u\in [0,1]$ and $j=1,2,\ldots,q+1$, we consider a probability $\zeta_j(u)$ on the nonnegative integers and for $1\leq j\leq q$, we denote by $\alpha_j(u)$ the mean of the distribution $\zeta_j(u)$.
Now let 
$$X_{n,k}=\sum_{j=1}^q \sum_{i=1}^{X_{n,k-j}}Y^{(n,k)}_{j,i}+\eta^{(n)}_k,$$
where for each integer $n\geq 1$, the family $\left\{Y^{(n,k)}_{j,i},\eta^{(n)}_h: (k,j,i,h)\in\Z^4\right\}$ contains independent 
random variables and such that for $1\leq j\leq q$, $Y^{(n,k)}_{j,i}$ has probability distribution $\zeta_j(k/n)$ and 
$\eta^{(n)}_k$ has probability distribution $\zeta_{q+1}(k/n)$. 
Note that, one can define a corresponding stationary autoregressive process. To this end, we denote by 
$F_{j,u}$ the cumulative distribution of the probability $\zeta_j(u)$ and we consider a family 
$\left\{U^{(k)}_{j,i},V_k: (i,j,k)\in\Z^3\right\}$ of i.i.d random variables uniformly distributed over $[0,1]$ and we set $Y^{(n,k)}_{j,i}=F_{j,k/n}^{-1}\left(U^{(k)}_{j,i}\right)$ where for a cumulative distribution function $G$, $G^{-1}$ denotes 
its left continuous inverse. Then one can consider the stationary version 
$$X_k(u)=\sum_{j=1}^q \sum_{i=1}^{X_{k-j}(u)} F_{j,u}^{-1}\left(U^{(k)}_{j,i}\right)+F_{q+1,u}^{-1}\left(V_k\right).$$ 

The following result is a consequence of Corollary \ref{autoreg}. 
Only the case $p=1$ is considered here.
\begin{cor}\label{INAR}
Let $d$ be the usual distance on $\R$.
Assume that $\sup_{u\in [0,1]}\sum_{j=1}^q\alpha_j(u)<1$, $\zeta_{q+1}$ has a finite first moment  and $W_1\left(\zeta_j(u),\zeta_j(v)\right)\leq C \vert u-v\vert^{\kappa}$ for $C>0,\kappa\in (0,1)$ and $1\leq j\leq q+1$.
Then Theorem \ref{central} and Proposition \ref{mixing2} apply.  
\end{cor}
\paragraph{Example.} Stationary INAR processes are often used when $\zeta_j$ is a Bernoulli distribution of parameter $\alpha_j\in (0,1)$ for $1\leq j\leq q$ and $\zeta_{q+1}$ is a Poisson distribution of parameter $\lambda$. This property guarantees 
that the marginal distribution is also Poissonian. Condition $\sum_{j=1}^q \alpha_j<1$ is a classical condition ensuring the existence of a stationary solution for this model.
In the locally stationnary case, let $U$ be a random variable following a uniform distribution over $[0,1]$ and $(N_t)_t$ be a Poisson process of intensity $1$. Then, we have 
$$W_1\left(\zeta_j(u),\zeta_j(v)\right)\leq \E\left\vert \mathds{1}_{U\leq \alpha_j(u)}-\mathds{1}_{U\leq \alpha_j(v)}\right\vert\leq \left\vert \alpha_j(u)-\alpha_j(v)\right\vert, $$
$$W_1\left(\zeta_{q+1}(u),\zeta_{q+1}(v)\right)\leq \E\left\vert N_{\lambda(u)}-N_{\lambda(v)}\right\vert\leq \left\vert \lambda(u)-\lambda(v)\right\vert.$$ 
Then the assumptions of Corollary \ref{INAR} are satisfied if the functions $\alpha_j$ and $\lambda$ are $\kappa-$Hölder continuous and if $\sum_{j=1}^q \alpha_j(u)<1$, $u\in [0,1]$. 

\paragraph{Note.} One can also state a result for $p\geq 1$. This case is important if we have to compare the expectation of some polynomials of the time-varying process with its the stationary version. However, in the example given above, a naive application of our results will require $\sum_{j=1}^q\alpha_j(u)^{1/p}<1$. Moreover, one can show that a $\kappa-$Hölder regularity on $\alpha_j$ and $\lambda$ entails 
a $\frac{\kappa}{p}-$Hölder regularity in Wasserstein metrics. 
For instance if $q=1$, we have $W_p\left(\delta_0Q_u,\delta_0Q_v\right)\geq \left\vert\lambda(u)-\lambda(v)\right\vert^{1/p}$. 
In order to avoid these unnatural conditions for this model, we will use the approach developed in Section \ref{Dob3}.

\section{Local stationarity and drift conditions}\label{Dob3}
In this section, we will use some drift and minoration conditions to extend the Dobrushin's contraction technique 
of Section \ref{Dob1}. A key result for this section is Lemma \ref{Mouli} which is adapted from 
Lemma $6.29$ in \citet{DM}. This result gives sufficient conditions for contracting Markov kernels 
with respect to norm induced by a particular Foster-Lyapunov drift function. The original argument for such contraction properties is due to \citet{HM1}. 
This important result will enable us to consider additional examples of locally stationary Markov chains with non compact state spaces. 
For a function $V:E\rightarrow [1,\infty)$, we define the $V-$norm of signed measure $\mu$ on $\left(E,\mathcal{B}(E)\right)$ 
by 
\begin{equation}\label{norm}
\Vert \mu\Vert_{V}=\sup_{f:\vert f(x)\vert\leq V(x)}\left\vert \int f d\mu\right\vert. 
\end{equation}

\subsection{General result}
We will assume that there exists a measurable function $V:E\rightarrow [1,\infty)$ such that
\begin{description}
\item[F1]
there exist $\epsilon>0$, $\lambda\in(0,1)$, an integer $m\geq 1$ and two real numbers $b>0,K\geq 1$ such that for all $u,u_1,\ldots,u_m\in [0,1]$ satisfying $\vert u-u_i\vert\leq \epsilon$, 
$$Q_u V\leq K V,\quad Q_{u_1}\cdots Q_{u_m} V\leq \lambda V+b.$$
Moreover, there exists $\eta>0$, $R>2b/(1-\lambda)$ and a probability measure $\nu\in\mathcal{P}(E)$ such that
$$\delta_x Q_{u_1}\cdots Q_{u_m}\geq \eta \nu,\mbox{ if } V(x)\leq R,$$

\item[F2]
there exist $\kappa\in (0,1)$ and a function $\widetilde{V}:E\rightarrow (0,\infty)$ such that $\sup_{u\in [0,1]}\pi_u \widetilde{V}<\infty$ and for all $x\in E$, $\Vert\delta_x Q_u-\delta_x Q_v\Vert_V\leq \widetilde{V}(x)\vert u-v\vert^{\kappa}$.
\end{description}

We first give some properties of the Markov kernels $Q_u$ with respect to the $V-$norm.
\begin{prop}\label{ergod}
Assume that assumptions ${\bf F1-F2}$ hold true.
\begin{enumerate}
\item
There exist $C>0$ and $\rho\in (0,1)$ such that for all $x\in E$, 
$$\sup_{u\in [0,1]}\Vert \delta_x Q_u^j-\pi_u\Vert_V\leq C V(x)\rho^j.$$ 
Moreover $\sup_{u\in [0,1]}\pi_u V<\infty$. 
\item
There exists $C>0$ such that for all $(u,v)\in [0,1]^2$,
$$\Vert\pi_u-\pi_v\Vert_V\leq C\left\vert u-v\right\vert^{\kappa}.$$
\end{enumerate}
\end{prop}

\paragraph{Proof of Proposition \ref{ergod}.} 
\begin{enumerate}
\item
According to Lemma \ref{Mouli}, there exists $(\gamma,\delta)\in (0,1)^2$ only depending $\lambda,b,\eta$, such that 
$$\Delta_{V_{\delta}}(Q_u^m)=\sup\left\{\frac{\Vert \mu R-\nu R\Vert_{V_{\delta}}}{\Vert \mu-\nu\Vert_{V_{\delta}}}:\quad \mu, \nu\in\mathcal{P}(E),\quad \mu V_{\delta}<\infty,\nu V_{\delta}<\infty\right\}\leq \gamma,$$
with $V_{\delta}=1-\delta+\delta V$.
From Theorem $6.19$ in \citet{DM} and Assumption ${\bf F1}$, we have a unique invariant probability for $Q_u$, satisfying $\pi_u V<\infty$ and for  $\mu\in\mathcal{P}(E)$ such that $\mu V<\infty$, we have
$$\Vert \mu Q_u^j-\pi_u\Vert_{V_{\delta}}\leq \max_{0\leq s\leq m-1}\Delta_{V_{\delta}}(Q_u^s)\gamma^{[j/m]}\Vert \mu-\pi_u\Vert_{V_{\delta}}.$$
Note that $\Vert \cdot\Vert_{V_{\delta}}\leq \Vert \cdot \Vert_V\leq \frac{1}{\delta}\Vert \cdot\Vert_{V_{\delta}}$ and the two norms are equivalent.
Using Lemma $6.18$ in \citet{DM}, we have 
$$\Delta_{V_{\delta}}(Q_u^s)=\sup_{x\neq y}\frac{\Vert \delta_x Q_u^s-\delta_y Q_u^s\Vert_{V_{\delta}}}{V_{\delta}(x)+V_{\delta}(y)}\leq \frac{K^s}{\delta}.$$
Then it remains to show that $\sup_{u\in [0,1]}\pi_u V<\infty$ or equivalently $\sup_{u\in [0,1]}\pi_u V_{\delta}<\infty$. 
But this a consequence of the contraction property of the application $\mu\mapsto \mu Q_u^m$ on the space $$\mathcal{M}_{\delta}=\left\{\mu\in\mathcal{P}(E): \mu V_{\delta}<\infty\right\}$$
endowed with the distance $d_{\delta}(\mu,\nu)=\Vert \mu-\nu\Vert_{V_{\delta}}$, which is a complete metric space (see Proposition $6.16$ in \citet{DM}). Hence we have 
$$\mu-\pi_u=\sum_{j=0}^{\infty}\left[\mu Q_u^{mj}-\mu Q_u^{m(j+1)}\right]$$
which defines a normally convergent series in $\mathcal{M}_{\delta}$ and 
$$\Vert \mu-\pi_u\Vert_{V_{\delta}}\leq \sum_{j=0}^{\infty}\gamma^j\Vert \mu-\mu Q_u^m\Vert_{V_{\delta}}\leq \frac{\mu V+K^m \mu V}{1-\gamma}.$$
This shows that $\sup_{u\in [0,1]}\pi_u V<\infty$ and 
the proof of the first point is now complete.

\item
To prove the second point, we decompose as in the two previous sections $\pi_u-\pi_v=\pi_uQ_u^m-\pi_v Q_u^m+\pi_v Q_u^m-\pi_v Q_v^m$.
This leads to the inequality 
$$\Vert \pi_u-\pi_v\Vert_{V_{\delta}}\leq \frac{\Vert \pi_vQ_u^m-\pi_v Q_v^m\Vert_{V_{\delta}}}{1-\gamma}.$$
Moreover, we have 
\begin{eqnarray*}
\Vert \pi_vQ_u^m-\pi_v Q_v^m\Vert_{V_{\delta}}&\leq& \Vert \pi_vQ_u^m-\pi_v Q_v^m\Vert_V\\
&\leq &\sum_{j=0}^{m-1}\Vert \pi_v Q_v^{m-j-1}(Q_v-Q_u)Q_u^j\Vert_V\\
&\leq& \sum_{j=0}^{m-1}K^j \Vert \pi_v Q_v^{m-j-1}(Q_v-Q_u)\Vert_V\\
&\leq& \sum_{j=0}^{m-1}K^j \cdot \pi_v \widetilde{V}\cdot\vert u-v\vert^{\kappa}.\\
\end{eqnarray*}
Hence the result follows with $C=\frac{\sup_{u\in[0,1]}\pi_u\widetilde{V}}{1-\gamma}\sum_{j=0}^{m-1}K^j$.$\square$  
\end{enumerate}

Now, we give our result about local stationarity.
\begin{theo}\label{FL1}
\begin{enumerate}
\item
Under the assumptions ${\bf F1-F2}$,  
there exists a positive real number $C$, only depending on $m,\lambda,b,K$ and $\sup_{u\in [0,1]}\pi_u \widetilde{V}$ such that 
$$\Vert\pi^{(n)}_k-\pi_u\Vert_V\leq C\left[\left\vert u-\frac{k}{n}\right\vert^{\kappa}+\frac{1}{n^{\kappa}}\right].$$
\item
Assume that assumptions ${\bf F1-F3}$ hold true and that in addition, for all $(u,v)\in [0,1]^2$, 
$$\Vert \delta_x Q_u-\delta_x Q_v\Vert_{TV}\leq L(x)\vert u-v\vert^{\kappa}\mbox{ with } G=\sup_{u\in [0,1]\atop 1\leq \ell'\leq\ell}\E L(X_{\ell}(u))V(X_{\ell'}(u)<\infty.$$ On $\mathcal{P}(E^j)$, we define 
the $V-$norm by
$$\Vert f\Vert_V=\sup\left\{\int fd\mu:\quad \vert f(x_1,\ldots,x_j)\vert\leq V(x_1)+\cdots+V(x_j)\right\}.$$
Then there exists $C_j>0$, not depending on $k,n,u$ and such that,
$$\Vert \pi^{(n)}_{k,j}-\pi_{u,j}\Vert_V\leq C_j\left[\left\vert u-\frac{k}{n}\right\vert^{\kappa}+\frac{1}{n^{\kappa}}\right].$$
\end{enumerate}
Moreover, the triangular array of Markov chains $\left\{X_{n,k}: n\in\Z^+, k\leq n\right\}$ is locally stationary. 
\end{theo}

\paragraph{Proof of Theorem \ref{FL1}.}
We assume that $n\geq \frac{m}{\epsilon}$.
\begin{enumerate}
\item
We start with the case $j=1$. Under the assumptions of the theorem, Lemma \ref{Mouli} guarantees the existence of $(\gamma,\delta)\in (0,1)^2$ such that for all $k\leq n$, $\Delta_{V_{\delta}}\left(Q_{\frac{k-m+1}{n}}\cdots Q_{\frac{k}{n}}\right)\leq \gamma$ with $V_{\delta}=1-\delta+\delta V$. In the sequel, we set $Q_{k,m}=Q_{\frac{k-m+1}{n}}\cdots Q_{\frac{k}{n}}$. Then we get 
\begin{eqnarray*}
\Vert \pi_k^{(n)}-\pi_u\Vert_{V_{\delta}}\leq \gamma\Vert\pi_{k-m}^{(n)}-\pi_u\Vert_{V_{\delta}}+\Vert \pi_u Q_{k,m}-\pi_u P_u^m\Vert_{V_{\delta}}.
\end{eqnarray*}
Now using ${\bf F1}$ and ${\bf F3}$, we get 
\begin{eqnarray*}
\Vert \pi_u Q_{k,m}-\pi_u Q_u^m\Vert_V&\leq& \sum_{j=0}^{m-1}\Vert \pi_u Q_u^{m-j-1} \left[Q_{\frac{k-j}{n}}-Q_u\right]Q_{\frac{k-j+1}{n}}\cdots Q_{\frac{k}{n}}\Vert_V\\
&\leq& \sum_{j=0}^{m-1}K^j\cdot \pi_u Q_u^{m-j+1}\widetilde{V}\cdot \left\vert u-\frac{k-j}{n}\right\vert^{\kappa}\\
&\leq& K^m \sup_{u\in [0,1]}\pi_u \widetilde{V} \sum_{s=k-m+1}^k\left\vert u-\frac{s}{n}\right\vert^{\kappa}.
\end{eqnarray*}
This shows the result in this case, by noticing that $\Vert\cdot\Vert_{V_{\delta}}\leq \Vert\cdot \Vert_{V}\leq \delta^{-1}\Vert\cdot\Vert_{V_{\delta}}$.

Now, if $n<m/\epsilon$, we have 
$$\Vert \pi_k^{n}-\pi_u\Vert_V\leq \pi_k^{(n)}V+\sup_{u\in[0,1]}\pi_u V\leq K^{\frac{m}{\epsilon}}\left(1+\sup_{u\in[0,1]}\pi_u V\right) \leq K^{\frac{m}{\epsilon}}\left(1+\sup_{u\in[0,1]}\pi_u V\right)\frac{m^{\kappa}}{\epsilon^{\kappa}}n^{-\kappa},$$
which leads to the result.

\item
Assume that the result is true for an integer $j\geq 1$. Let $f:E^{j+1}\rightarrow \R_+$ be such that 
$f(x_1,\ldots,x_{j+1})\leq V(x_1)+\cdots+V(x_{j+1})$. Setting $s_j=\sum_{i=1}^j V(x_i)$ and  $g_j(x_{j+1})=f(x_1,\ldots,x_{j+1})$, we use the decomposition 
$$g_j=g_j\mathds{1}_{g_j\leq  V}+\left(g_j-V\right)\mathds{1}_{V< g_j\leq V+s_j}+V\mathds{1}_{V< g_j\leq V+s_j},$$
and we get 
\begin{eqnarray*}
\left\vert \delta_{x_{j+1}} Q_{\frac{k+j+1}{n}}g_j-\delta_{x_{j+1}} Q_u g_j\right\vert&\leq& \left\Vert \delta_{x_{j+1}} Q_{\frac{k+j+1}{n}}-\delta_{x_{j+1}} Q_u \right\Vert_V+ s_j\Vert\delta_{x_{j+1}} Q_{\frac{k+j+1}{n}}-\delta_{x_{j+1}} Q_u\Vert_{TV}\\
&\leq& \left[\widetilde{V}\left(x_{j+1}\right)+\left(V(x_1)+\cdots+V(x_j)\right)L(x_j)\right]\left\vert u-\frac{k+j+1}{n}\right\vert^{\kappa}.
\end{eqnarray*}
This yields to 
$$\Vert \pi_{u,j}\otimes Q_{\frac{k+j+1}{n}}-\pi_{u,j+1}\Vert_V\leq 2\left(\sup_{u\in [0,1]}\pi_u\widetilde{V}+j G\right)\left\vert u-\frac{k+j+1}{n}\right\vert^{\kappa}.$$
On the other hand 
$$\Vert \pi_{u,j}\otimes Q_{\frac{k+j+1}{n}}-\pi^{(n)}_{k,j+1}\Vert_V\leq (1+K)\Vert \pi^{(n)}_{k,j}-\pi_{u,j}\Vert_V.$$
The two last bounds lead to the result using finite induction. Moreover, using the same type of arguments, one can check the continuity condition $1$ of Definition \ref{debase}.$\square$
\end{enumerate}

One can also define a useful upper bound of the usual $\beta-$mixing coefficient which is useful to control covariances of unbounded functionals of the Markov chain. More precisely, we set 
$$\beta^{(V)}_n(j)= \sup_{k\leq n}\E\Vert \pi_k^{(n)}- \delta_{X_{n,k-j}}Q_{\frac{k-j+1}{n}}\cdots Q_{\frac{k}{n}}\Vert_V.$$
We have the following result which proof is straigthforward.
\begin{prop}\label{mixsuf}
Assume that assumption ${\bf F1}$ holds true and that $n\geq m/\epsilon$.
Then if $j=mg+s$, we have $\beta^{(V)}_n(j)\leq 2\delta^{-1}\sup_{k\leq n}\pi_k^{(n)}V\cdot K^s \gamma^g$,
where $\delta,\gamma\in (0,1)$ are given in Lemma \ref{Mouli}.   
\end{prop}
\paragraph{Notes}
\begin{enumerate}
\item
From the drift condition in ${\bf F1}$, we have $\pi_k^{(n)}V\leq \frac{b}{1-\lambda}$ for $n\geq m/\epsilon$. Hence $\sup_{n\in\Z^+, k\leq n}\pi^{(n)}_k V<\infty$. 
\item
We did not adapt the notion of $V-$mixing given in \citet{MT}, Chapter $16$. However, let us mention that if $A=\sup_{k\leq n,n\in\Z^+}\pi^{(n)}_k \left[\vert g\vert V\right]<\infty$, we get the following covariance inequality
$$\left\vert \c\left(g\left(X_{n,k-j}\right),f\left(X_{n,k}\right)\right)\right\vert\leq 2\delta^{-1}A K^s \gamma^g,$$
when $j=mg+s$.
\end{enumerate}
\subsection{Example $1$: the random walk on the nonnegative integers}
Let $p,q,r:[0,1]\rightarrow(0,1)$ three $\kappa-$Hölder continuous functions such that $p(u)+q(u)+r(u)=1$ and $\frac{p(u)}{q(u)}<1$.
For $x\in \N^*$, we set $Q_u(x,x)=r(u)$, $Q_u(x,x+1)=p(u)$ and $Q_u(x,x-1)=q(u)$. Finally 
$Q_u(0,1)=1-Q_u(0,0)=p(u)$.  
In the homogeneous case, geometric ergodicity holds under the condition $p<q$. See \citet{MT}, 
Chapter $15$. In this case the function $V$ defined by $V(x)=z^x$ is a Foster-Lyapunov function if $1<z<q/p$. For the non-homogeneous case, let $z\in (1,e)$ where $e=\min_{u\in [0,1]}q(u)/p(u)$. We set $\gamma=\max_{u\in [0,1]}\left\{r(u)+p(u)z+q(u)z^{-1}\right\}$ and $\overline{p}=\max_{u\in [0,1]}p(u)$.
Note that 
$$\gamma\leq 1+\overline{p}(z-1)\max_{u\in [0,1]}\left[1-\frac{q(u)}{p(u)z}\right]\leq 1+\overline{p}(z-1)\left[1-\frac{e}{z}\right]<1.$$
Then we have $Q_uV(x)\leq \gamma V(x)$ for all $x>0$ and $Q_uV(0)=p(u)z+(1-p(u))\leq c=\overline{p}(z-1)+1$. 
For an integer $m\geq 1$, we have $Q_{u_1}\cdots Q_{u_m} V\leq \gamma^m V+\frac{c}{1-\gamma}$. If $m$ is large enough, we have 
$\frac{2c}{(1-\gamma)(1-\gamma^m)V(m)}<1$. Moreover, for such $m$, if $R=V(m)$, we have $\{V\leq R\}=\{0,1,\ldots,m\}$ and if $x=0,\ldots,m$, we have 
$\delta_xQ_{u_1}\cdots Q_{u_m}\geq \eta \delta_0$ for a $\eta>0$. Assumption ${\bf F3}$ is immediate. Moreover the additional condition in the second point of Theorem \ref{FL1} is automatically checked with a constant function $L$.\\
However this example is more illustrative. Indeed parameters $p(u)$ and $q(u)$ can be directly estimated by 
$$\hat{p}(u)=\sum_{i=1}^{n-1}e_i(u)\mathds{1}_{X_{n,i+1}-X_{n,i}=1},\quad \hat{q}(u)=\sum_{i=1}^{n-1}e_i(u)\mathds{1}_{X_{n,i+1}-X_{n,i}=-1},$$
where the weights $e_i(u)$ are defined as in Subsection \ref{poid}. The indicators are independent Bernoulli random variables with parameter $p\left(\frac{i+1}{n}\right)$ or $q\left(\frac{i+1}{n}\right)$ and the asymptotic behavior of the estimates is straightforward.

\subsection{Example $2$: INAR processes}
We again consider INAR processes. For simplicity, we only consider the case $q=1$ with Bernoulli counting sequences and a Poissonian noise. The parameters $\alpha(u)$ (resp. $\lambda(u)$) of the counting sequence (resp. the Poissonian noise) are assumed to be $\kappa-$Hölder continuous. We will show that our results apply with drift functions $V_p(x)=x^p+1$ for an arbitrary integer $p\geq1$. To this end, we consider a sequence $\left(Y_i(u)\right)_{i\geq 0}$ of i.i.d random variables following the Bernoulli distribution of parameter $\alpha(u)$ and a random variable $\xi(u)$ following the Poisson distribution of parameter $\lambda(u)$. We assume that $\xi(u)$ and the sequence $\left(Y_i(u)\right)_{i\geq 0}$ are independent. For $u\in [0,1]$, we have 
\begin{eqnarray*}
\delta_x Q_u V_p&=& 1+\E\left(\alpha(u)x+\sum_{i=1}^x\left(Y_i(u)-\alpha(u)\right)+\xi(u)\right)^p\\
&=& 1+\alpha(u)^p x^p+\sum_{j=1}^p\begin{pmatrix} p\\j\end{pmatrix}\alpha(u)^{p-j}x^{p-j}\E\left(\sum_{i=1}^x \left(Y_i(u)-\alpha(u)\right)+\xi(u)\right)^j\\
\end{eqnarray*}
Using the Burkhölder inequality for martingales, we have for an integer $\ell\geq 2$, 
$$\E\left(\sum_{i=1}^x\left(Y_i(u)-\alpha(u)\right)\right)^{\ell}\leq C_{\ell}x^{\frac{\ell}{2}}\max_{1\leq i\leq x}\E\left\vert Y_i(u)-\alpha(u)\right\vert^{\ell}\leq C_{\ell}x^{\frac{\ell}{2}},$$
where $C_{\ell}$ is a universal constant. Then, we deduce from the previous equalities that there exist two constants $N_1$ and $N_2$ such that 
$$\delta_x Q_u V_p\leq \alpha(u)^p V_p(x)+M_1 x^{p-1}+M_2.$$
To check the drift condition in ${\bf F1}$ for $m=1$, one can choose $\gamma>0$ such that $\lambda=\max_{u\in [0,1]}\alpha(u)^p+\gamma<1$ and $b=M_2+\left(\frac{\gamma}{M_1}\right)^{p-1}$.
In this case, the minoration condition is satisfied on each finite set $\mathcal{C}$ with $\nu=\delta_0$ because 
$$Q_u(x,0)\geq \left(1-\max_{u\in [0,1]}\alpha(u)\right)^x\exp\left(-\max_{u\in [0,1]}\lambda(u)\right)$$ and $\eta=\min_{u\in [0,1]}\min_{x\in\mathcal{C}}Q_u(x,0)>0$.
This shows that assumption ${\bf F1}$ is satisfied by taking $R$ large enough.\\
Finally, we show ${\bf A3}$. Let $u,v\in [0,1]$. Denoting $\overline{\lambda}=\max_{u\in [0,1]}\lambda(u)$ and by $\mu_u$ the Poisson distribution of parameter $\lambda(u)$, we have
$$\max_{u\in [0,1]}\mu_u V_p\leq 1+\E N_{\overline{\lambda}}^p,\quad \Vert \mu_u-\mu_v\Vert_{V_p}\leq \sum_{k\geq 0}\frac{V_p(k)}{k!}\left(k\bar{\lambda}^{k-1}+\bar{\lambda}^k\right)\cdot\left\vert \lambda(u)-\lambda(v)\right\vert,$$
where $(N_t)_{t\geq 0}$ is Poisson process of intensity $1$. Moreover, if $\nu_u$ denotes the Bernoulli distribution of parameter $\alpha(u)$, we have $\Vert \nu_u-\nu_v\Vert_{V_p}\leq 3\left\vert \alpha(u)-\alpha(v)\right\vert$. From Lemma \ref{interm}, we easily deduce that ${\bf F2}$ holds for $\widetilde{V}=C V_{p+1}$ where $C$ is a positive real number. Note that we have $\sup_{u\in [0,1]}\pi_u \widetilde{V}<\infty$ because $\widetilde{V}$ also satisfies the drift and minoration condition.    

Let us now give an estimator for parameter $\left(\alpha(u),\lambda(u)\right)$. 
A natural estimate is obtained by localized least squares. Setting $a(u)=\left(\alpha(u),\lambda(u)\right)'$ and 
$\mathcal{Y}_{n,i}=\left(1,X_{n,i-1}\right)'$. Then we define 
$$\hat{a}(u)=\arg\min_{\alpha}\sum_{i=2}^n e_i(u)\left(X_{n,i}-\mathcal{Y}_{n,i}'\alpha\right)^2=\left(\sum_{i=2}^n e_i(u)\mathcal{Y}_{n,i}\mathcal{Y}_{n,i}'\right)^{-1}\sum_{i=2}^ne_i(u)X_{n,i}\mathcal{Y}_{n,i},$$
where the weights $e_i(u)$ were defined in Subsection \ref{poid}.
Using our results and assuming that $b\rightarrow 0$ and $nb\rightarrow \infty$, we get
$$\sum_{i=2}^n e_i(u)\E\mathcal{Y}_{n,i}\mathcal{Y}_{n,i}'=\E\mathcal{Y}_i(u)\mathcal{Y}_i(u)'+O\left(b^{\kappa}+\frac{1}{n^\kappa}\right)=O\left(b^{\kappa}\right).$$
In the same way, we have 
$$\sum_{i=2}^n e_i(u)\E X_{n,i}\mathcal{Y}_{n,i}=\E X_i(u)\mathcal{Y}_i(u)+O\left(b^{\kappa}+\frac{1}{n^\kappa}\right)=O\left(b^{\kappa}\right).$$
Moreover, using our covariance inequality (see the notes after Proposition \ref{mixsuf}), we get 
$$\v\left(\sum_{i=2}^n e_i(u)\mathcal{Y}_{n,i}\mathcal{Y}_{n,i}'\right)=O\left((nb)^{-1}\right).$$
Moreover using the decomposition $X_{n,i}=\mathcal{Y}_{n,i}'a(i/n)+X_{n,i}-\E\left(X_{n,i}\vert \mathcal{F}_{n,i-1}\right)$ where $\mathcal{F}_{n,i}=\sigma\left(X_{n,j}:j\leq i\right)$ and the fact that for all $p\geq 1$, $\sup_{n\in\Z^+,k\leq n}\E\vert X_{n,k}\vert^p<\infty$, we also obtain 
$$\v\left(\sum_{i=2}^n e_i(u)\mathcal{Y}_{n,i}X_{n,i}\right)=O\left((nb)^{-1}\right).$$
Collecting all the previous properties, we get $\hat{a}(u)=a(u)+O_{\P}\left(b^{\kappa}+\frac{1}{\sqrt{nb}}\right)$. Asymptotic normality or uniform control of $\hat{a}(u)-a(u)$ can also be obtained using adapted results for strong mixing sequences.

\section{Auxiliary results}

\subsection{Auxiliary result for Section \ref{Dob1}}
\begin{prop}\label{TLC}
Let $\left(H_i^{(n)}\right)_{1\leq i\leq n,n>0}$ be a double array of real-valued random variables with finite variance and mean zero.
Let $\left(\alpha_n(k)\right)_{k\geq 0}$ be the sequence of strong mixing coefficients of the sequence $\left(H_i^{(n)}\right)_{1\leq i\leq n,n>0}$ and $\alpha_{(n)}^{-1}$ be the inverse function of the associated mixing rate function. Suppose that 
\begin{equation}\label{condvar}
\lim\sup_{n\rightarrow \infty}\max_{1\leq i\leq n}\frac{V_{n,i}}{V_{n,n}}<\infty,
\end{equation}
where $V_{n,i}=\v\left(\sum_{j=1}^iH^{(n)}_j\right)$. Let 
$$Q_{n,i}=\sup\left\{t\in \R_+: \P\left(\left\vert H_i^{(n)}\right\vert>t\right)>u\right\}.$$
Then $\sum_{i=1}^n H_i^{(n)}$ converges to the standard normal distribution if 
\begin{equation}\label{condmix}
V_{n,n}^{-3/2}\sum_{i=1}^n\int_0^1 \alpha_{(n)}^{-1}(x/2)Q_{n,i}^2(x)\inf\left(\alpha_{(n)}^{-1}(x/2)Q_{n,i}(x),\sqrt{V_{n,n}}\right)dx\rightarrow 0,
\end{equation}
as $n$ tends to $\infty$.
\end{prop}

\subsection{Auxiliary Lemmas for Section \ref{Dob2}}
\begin{lem}\label{preli}
Let $\mu\in\mathcal{P}_p(E)$ and $Q$, $R$ be two probability kernels from $\left(E,\mathcal{B}(E)\right)$ to $\left(E,\mathcal{B}(E)\right)$  such that
\begin{enumerate}
\item
for all $x\in E$,  the two probability measures $\delta_x Q$ and $\delta_x R$ are  elements of $\mathcal{P}_p(E)$,
\item
there exists $C>0$ such that for all $(x,y)\in E^2$, 
$$W_p\left(\delta_x Q,\delta_y Q\right)\leq C d(x,y),\quad W_p\left(\delta_x R,\delta_y R\right)\leq C d(x,y).$$ 
\end{enumerate} 
Then, if $\mu\in\mathcal{P}_p(E)$, the two probability measures $\mu Q$, $\mu R$ are also elements of $\mathcal{P}_p(E)$. Moreover, we have 
\begin{equation}\label{first1}
W^p_p\left(\mu Q,\mu R\right)\leq \int W^p_p\left(\delta_x Q,\delta_x R\right)d\mu(x),
\end{equation}
and if $\nu$ is another element of $\mathcal{P}_p(E)$, we have 
\begin{equation}\label{second1}
W_p\left(\mu Q,\nu Q\right)\leq C W_p(\mu,\nu).
\end{equation}
\end{lem}

\paragraph{Proof of Lemma \ref{preli}.}
Using Lemma \ref{Lip} for $f(x)=d(x,x_0)$, we have for a given $y\in E$,
\begin{eqnarray*}
\int d(x,x_0)^p Q(y,dx)&\leq& \left[W_p(\delta_y Q,\delta_{x_0} Q)+\left(\int d(x,x_0)^p Q(x_0,dx)\right)^{1/p}\right]^p\\
&\leq& \left[C d(x_0,y)+\left(\int d(x,x_0)^p Q(x_0,dx)\right)^{1/p}\right]^p.\\
\end{eqnarray*}
After integration with respect to $\mu$, it is easily seen that $\mu Q\in\mathcal{P}_p(E)$.\\
To show (\ref{first1}), one can use Kantorovitch duality (see \citet{Villani}, Theorem $5.10$). Denoting by $\mathcal{C}_b(E)$ the set of bounded continuous functions on $E$, we have
\begin{eqnarray*} 
W^p_p\left(\mu Q,\mu R\right)&=&\sup_{\phi(x)-\psi(y)\leq d(x,y)^p, (\phi,\psi)\in \mathcal{C}_b(E)}\left\{\int \phi(x)\mu Q(dx)-\int \psi(y)\mu R(dy)\right\}\\
&\leq& \int \left[\sup_{\phi(x)-\psi(y)\leq d(x,y)^p, (\phi,\psi)\in \mathcal{C}_b(E)}\left\{\int \phi(x)Q(z,dx)-\int \psi(y) R(z,dy)\right\}\right]\mu(dz)\\
&\leq & \int W^p_p\left(\delta_z Q,\delta_z R\right)\mu(dz).
\end{eqnarray*}
Finally, we show (\ref{second1}). Let $\phi,\psi$ be two elements of $\mathcal{C}_b(E)$ such that $\phi(x)-\psi(y)\leq d(x,y)^p$ and $\gamma$ an optimal coupling for $(\mu,\nu)$. 
Then, for $u,v\in E$, we have  
$$\int \phi(x) Q(u,dx)-\int\psi(y)Q(v,dy)\leq W^p_p\left(\delta_u Q,\delta_v Q\right)\leq C^p d(u,v)^p.$$
Moreover, 
$$\int \phi(x)\mu Q(dx)-\int \psi(y)\nu Q(dy)=\int \gamma(du,dv)\left[\int \phi(x)Q(u,dx)-\int \psi(y)Q(v,dy)\right].$$ 
Then (\ref{second1}) easily follows from Kantorovitch duality.$\square$

\begin{lem}\label{preli2}
Let $j\geq 1$ be an integer. Assume that $Q_1,\ldots,Q_j$  and $R_1,\ldots,R_j$ are Markov kernels such that for all $x\in E$ and $1\leq i\leq j$, $\delta_x Q_i$ and $\delta_x R_i$ are elements of $\mathcal{P}_p(E)$ satisfying 
$$W_p\left(\delta_x Q_i,\delta_y Q_i\right)\leq L_id(x,y),\quad W_p\left(\delta_x R_i,\delta_y R_i\right)\leq L_id(x,y),$$
for all $(x,y)\in E^2$. Then, for all $x\in E$, we have
$$W_p\left(\delta_x Q_1\cdots Q_j,\delta_x R_1\cdots R_j\right)\leq \sum_{s=0}^{j-1}L_j\cdots L_{j-s+1}D_{j-s},$$
where $D^p_i=\int W_p^p\left(\delta_y Q_i,\delta_y R_i\right)\delta_x R_1\cdots R_{i-1}(dy)$.
\end{lem}
\paragraph{Proof of Lemma \ref{preli2}}
Using the inequality 
$$W_p\left(\delta_x Q_1\cdots Q_j,\delta_x R_1\cdots R_j\right)\leq \sum_{s=0}^{j-1}W_p\left(\delta_xR_1\cdots R_{j-s-1}Q_{j-s}\cdots Q_j,
\delta_xR_1\cdots R_{j-s}Q_{j-s+1}\cdots Q_j\right),$$
the result follows using Lemma \ref{preli}.$\square$

\begin{lem}\label{Lip}
If $f:E\rightarrow \R$ is a Lipschitz function, then for all measures $\mu,\nu\in \mathcal{P}_p(E)$, we have 
$$\left\vert \left(\int f^pd\mu\right)^{1/p}-\left(\int f^pd\nu\right)^{1/p}\right\vert\leq \delta(f)W_p(\mu,\nu),$$
where $\delta(f)$ denotes the Lipschitz constant of $f$:
$$\delta(f)=\sup_{x\neq y}\frac{\vert f(x)-f(y)\vert}{d(x,y)}.$$
\end{lem}

\paragraph{Proof of Lemma \ref{Lip}}

If $\gamma$ denotes an optimal coupling for $(\mu,\nu)$, we get from the triangular inequality,
\begin{eqnarray*}
&&\left\vert \left(\int f^pd\mu\right)^{1/p}-\left(\int f^pd\nu\right)^{1/p}\right\vert\\
&=& \left\vert \left(\int f^p(x)d\gamma(x,y)\right)^{1/p}-\left(\int f^p(y)d\gamma(x,y)\right)^{1/p}\right\vert\\
&\leq& \left(\int \vert f(x)-f(y)\vert^pd\gamma(x,y)\right)^{1/p}\\
&\leq& \delta(f)\left(\int d(x,y)^pd\gamma(x,y)\right)^{1/p}.
\end{eqnarray*}
which leads to the result of the lemma.$\square$

\begin{lem}\label{astuce}
Let $X$ and $Y$ two random variables taking values in $(E,d)$ and such that $\P_X$, $\P_Y\in \mathcal{P}_d(E)$. 
On $E\times E$, we define the metric 
$$\widetilde{d}\left((x_1,x_2),(y_1,y_2)\right)=\left(d(x_1,y_1)^p+d(x_2,y_2)^p\right)^{1/p}.$$
Then we have
$$W_p\left(\P_{X,Y},\P_{Y,Y}\right)\geq 2^{-\frac{p-1}{p}}\E^{1/p}\left(d(X,Y)^p\right).$$
\end{lem}

\paragraph{Proof of Lemma \ref{astuce}}
Consider the Lipschitz function $f:E\times E\rightarrow \R$ defined by $f(x_1,x_2)=d(x_1,x_2)$. Using the triangular inequality and convexity, we have $\delta(f)\leq 2^{\frac{p-1}{p}}$.
Then the result is a consequence of Lemma \ref{Lip}.$\square$

\subsection{Auxiliary Lemmas for Section \ref{Dob3}}
The following result is an adaptation of Lemma $6.29$ given in \citet{DM}. The proof is omitted because the arguments are exactly the same. See also \citet{HM1} for the original proof of this result.
Note however, that we use the condition $R>\frac{2b}{1-\lambda}$ instead of $R>\frac{2b}{(1-\lambda)^2}$ because coefficients $(\lambda,b)$ are obtained for $m$ iterations of the kernel (in \citet{DM}, these coefficients are that for the case $m=1$).
For a Markov kernel $R$ on $\left(E,\mathcal{B}(E)\right)$, we define its Dobrushin's contraction coefficient by 
$$\Delta_V(R)=\sup\left\{\frac{\Vert \mu R-\nu R\Vert_V}{\Vert \mu-\nu\Vert_V}:\quad \mu,\nu\in\mathcal{P}(E),\quad \mu V<\infty,\nu V<\infty\right\}.$$
\begin{lem}\label{Mouli}
Under assumption ${\bf F1}$, there exists $(\gamma,\delta)\in (0,1)^2$, only depending on $\lambda,\eta,b$, such that for all $(u,u_1,\ldots,u_m)\in [0,1]^{m+1}$ 
such that $\vert u_i-u\vert\leq \epsilon$, $1\leq i\leq m$, we have
$$\Delta_{V_{\delta}}\left(Q_{u_1}\cdots Q_{u_m}\right)\leq \gamma,\mbox{ where } V_{\delta}=1-\delta+\delta V.$$
\end{lem}

\begin{lem}\label{interm}
Let $X_1,X_2,\ldots,X_n,Y_1,Y_2,\ldots,Y_n$ be independent random variables such that $A_n=\max_{1\leq i\leq n}\E V(X_i)\vee\E V(Y_i)<\infty$ for $1\leq i\leq n$, with $V(x)=1+\vert x\vert^p$ and $p\geq 1$. Then we have 
$$\sup_{\vert f\vert \leq V}\left\vert \E f(X_1+\cdots+X_n)-\E f(Y_1+\cdots+Y_n)\right\vert\leq 2^{p+1}n^{p+1}\cdot A_n\cdot \max_{1\leq i\leq n}\Vert \P_{X_i}-\P_{Y_i}\Vert_V.$$
\end{lem}

\paragraph{Proof of Lemma \ref{interm}.}
Note first that if $\vert f(x)\vert\leq V(x)$ for all $x\in E$, then 
$\vert f(x+y)\vert\leq 2^pV(x)V(y)$. This leads to  
\begin{eqnarray*}
&&\left\vert \E f(X_1+\cdots+X_n)-\E f(Y_1+\cdots+Y_n)\right\vert\\
&\leq& \sum_{j=1}^n \left\vert \E f\left(X_1+\cdots+X_{j-1}+X_j+Y_{j+1}+\cdots+Y_n\right)-f\left(X_1+\cdots+X_{j-1}+Y_j+Y_{j+1}+\cdots+Y_n\right)\right\vert\\
&\leq& 2^p \sum_{j=1}^n \Vert \P_{X_j}-\P_{Y_j}\Vert_V\cdot \E V\left(X_1+\cdots+X_{j-1}+Y_{j+1}+\cdots+Y_n\right)\\
&\leq& 2^p (n-1)^p \sum_{j=1}^n \Vert \P_{X_j}-\P_{Y_j}\Vert_V A_n\\
&\le& 2^p n^{p+1}A_n \max_{1\leq i\leq n}\Vert \P_{X_i}-\P_{Y_i}\Vert_V.\square
\end{eqnarray*}

\bibliographystyle{plainnat}
\bibliography{bibMarkov}

\begin{thebibliography}{25}
\providecommand{\natexlab}[1]{#1}
\providecommand{\url}[1]{\texttt{#1}}
\expandafter\ifx\csname urlstyle\endcsname\relax
  \providecommand{\doi}[1]{doi: #1}\else
  \providecommand{\doi}{doi: \begingroup \urlstyle{rm}\Url}\fi

\bibitem[Al~Osh and Alzaid(1987)]{Al}
M.~Al~Osh and A.~Alzaid.
\newblock Firs-order integer-valued autoregressive process.
\newblock \emph{J. Time Series Anal.}, 8:\penalty0 261--275, 1987.

\bibitem[Dahlhaus(1997)]{Dahlhaus}
R.~Dahlhaus.
\newblock Fitting time series models to nonstationary processes.
\newblock \emph{Ann. Statist.}, 25:\penalty0 1--37, 1997.

\bibitem[Dahlhaus and Subba~Rao(2006)]{DR}
R.~Dahlhaus and S.~Subba~Rao.
\newblock Statistical inference for time-varying arch processes.
\newblock \emph{Ann. Statist.}, 34:\penalty0 1075--1114, 2006.

\bibitem[Dedecker and Prieur(2004)]{ded}
J.~Dedecker and C.~Prieur.
\newblock Coupling for $\tau-$dependent sequences and applications.
\newblock \emph{Journal of Theoretical Probability}, 17:\penalty0 861--885,
  2004.

\bibitem[Dobrushin(1956)]{Dob}
R.L. Dobrushin.
\newblock Central limit theorems for nonstationary markov chains.
\newblock \emph{Th. Prob. Appl.}, 1:\penalty0 329--383, 1956.

\bibitem[Dobrushin(1970)]{Dob+}
R.L. Dobrushin.
\newblock Prescribing a system of random variables by conditional
  distributions.
\newblock \emph{Th. Prob. Appl.}, 15:\penalty0 458--486, 1970.

\bibitem[Douc et~al.(2004)Douc, Moulines, and Rosenthal]{Mouli}
R.~Douc, E.~Moulines, and J.S. Rosenthal.
\newblock Quantitative bounds on convergence of time-inhomogeneous markov
  chains.
\newblock \emph{Ann. Appl. Probab.}, 14(4):\penalty0 1643--1665, 2004.

\bibitem[Douk et~al.(2014)Douk, Moulines, and Stoffer]{DM}
R.~Douk, E.~Moulines, and D.~Stoffer.
\newblock \emph{Nonlinear Time Series}.
\newblock Chapman and Hall, 2014.

\bibitem[Doukhan(1994)]{Do}
P.~Doukhan.
\newblock \emph{Mixing. Properties and Examples}.
\newblock Springer-Verlag, 1994.

\bibitem[Francq and Zako\"ian(2010)]{FZ}
C.~Francq and J-M. Zako\"ian.
\newblock \emph{GARCH models: structure, statistical inference and financial
  applications}.
\newblock Wiley, 2010.

\bibitem[Fryzlewicz et~al.(2008)Fryzlewicz, Sapatinas, and Subba~Rao]{Fryz}
P.~Fryzlewicz, T.~Sapatinas, and S.~Subba~Rao.
\newblock Normalized least-squares estimation in time-varying arch models.
\newblock \emph{Ann. Statist.}, 36:\penalty0 742--786, 2008.

\bibitem[Hairer and Mattingly(2011)]{HM1}
M.~Hairer and J.C. Mattingly.
\newblock Yet another look at harris' ergodic theorem for markov chains.
\newblock In \emph{Seminar on Stochastic Analysis, Random Fields and
  Applications IV}, volume~63, pages 109--117. Birkhäuser,/Springer Basel AG,
  Basel, 2011.

\bibitem[Jin-Guan and Yuan(1991)]{Guan}
D.~Jin-Guan and L.~Yuan.
\newblock The integer-valued autoregressive (inar(p)) model.
\newblock \emph{J. Time Series Anal.}, 12:\penalty0 129--142, 1991.

\bibitem[Meyn and Tweedie(2009)]{MT}
S.~Meyn and R.L. Tweedie.
\newblock \emph{Markov Chains and Stochastic Stability 2nd}.
\newblock Cambridge University Press New York, 2009.

\bibitem[Rio(1995)]{Rio2}
E.~Rio.
\newblock About the lindeberg method for strongly mixing sequences.
\newblock \emph{ESAIM, Probability and Statistics}, 1:\penalty0 35--61, 1995.

\bibitem[Rio(1999)]{Rio1}
E.~Rio.
\newblock \emph{Théorie asymptotique des processus aléatoires faiblement
  dépendants}.
\newblock Springer, 1999.

\bibitem[Rio(2013)]{Rio11}
E.~Rio.
\newblock Inequalities and limit theorems for weakly dependent sequences.
\newblock \url{https://cel.archives-ouvertes/cel-00867106}, 2013.

\bibitem[Saloff-Coste and Z\'u\~{n}iga(2007)]{SC1}
L.~Saloff-Coste and J.~Z\'u\~{n}iga.
\newblock Convergence of some time-inhomogeneous markov chains via spectral
  techniques.
\newblock \emph{Stochastic Process. Appl.}, 117:\penalty0 961--979, 2007.

\bibitem[Saloff-Coste and Z\'u\~{n}iga(2011)]{SC2}
L.~Saloff-Coste and J.~Z\'u\~{n}iga.
\newblock Merging for inhomogeneous finite markov chains, part ii: Nash and
  log-sobolev inequalities.
\newblock \emph{Ann. Probab.}, 39:\penalty0 1161--1203, 2011.

\bibitem[Subba~Rao(2006)]{SR}
S.~Subba~Rao.
\newblock On some nonstationary, nonlinear random processes and their
  stationary approximations.
\newblock \emph{Adv. in App. Probab.}, 38:\penalty0 1155--1172, 2006.

\bibitem[Truquet(2016)]{Truq}
L.~Truquet.
\newblock Parameter stability and semiparametric inference in time-varying arch
  models.
\newblock \emph{Forthcoming in JRSSB}, 2016.

\bibitem[Villani(2009)]{Villani}
C.~Villani.
\newblock \emph{Optimal Transport. Old and New}.
\newblock Springer, 2009.

\bibitem[Vogt(2012)]{Vogt}
M.~Vogt.
\newblock Nonparametric regression for locally stationary time series.
\newblock \emph{Ann. Statist.}, 40:\penalty0 2601--2633, 2012.

\bibitem[Winkler(1995)]{Winkler}
G.~Winkler.
\newblock \emph{Image Analysis, Random Fields and Dynamic Monte Carlo Methods}.
\newblock Springer, 1995.

\bibitem[Zhang and Wu(2015)]{ZW}
T.~Zhang and W.B. Wu.
\newblock Time-varying nonlinear regression models: nonparametric estimation
  and model selection.
\newblock \emph{Ann. Statist.}, 43:\penalty0 741--768, 2015.

\end{thebibliography}

\end{document}